\documentclass[11pt]{amsart}

\usepackage{latexsym,amsfonts,amssymb,exscale,enumerate,comment}
\usepackage{amsmath,amsthm,amscd}
\usepackage[all,knot]{xy}

\usepackage{fullpage, braket}

\usepackage{caption}
\usepackage{subcaption}

\usepackage{url}
\usepackage[bookmarks=true,%
    colorlinks=true,%
    linkcolor=blue,%
    citecolor=blue,%
    filecolor=blue,%
    menucolor=blue,%
    urlcolor=blue,%
    breaklinks=true]{hyperref}

\input xy
\usepackage[all]{xy}
\SelectTips{cm}{}

\usepackage{tikz}
\usepackage{tikz-cd} 
\usetikzlibrary{shapes}
\usetikzlibrary{decorations.markings}
\usetikzlibrary{decorations.pathreplacing}
\tikzstyle directed=[postaction={decorate,decoration={markings,
    mark=at position #1 with {\arrow{>}}}}]
\tikzstyle rdirected=[postaction={decorate,decoration={markings,
    mark=at position #1 with {\arrow{<}}}}]
    
\newcommand{\hackcenter}[1]{
 \xy (0,0)*{#1}; \endxy}

\tikzset{->-/.style={decoration={
  markings,
  mark=at position #1 with {\arrow{>}}},postaction={decorate}}}

\usepackage{graphicx}
\usepackage{color}
\newcommand{\brk}[1]{{\left\langle{#1}\right\rangle}}

\newcommand{\md}{\operatorname{\mathsf{d}}}

\newcommand{\SL}{\mathrm{SL}}

\newcommand{\cat}{\mathcal{C}}

\theoremstyle{plain}
\newtheorem{theorem}{Theorem}

\newtheorem*{theo}{Theorem}

\newtheorem{proposition}[theorem]{Proposition}
\newtheorem{lemma}[theorem]{Lemma}

\theoremstyle{definition}

\theoremstyle{definition}
\newtheorem{remark}[theorem]{Remark}

\newcommand{\maps}{\colon}

\newcommand{\refequal}[1]{\xy {\ar@{=}^{#1}
(-1,0)*{};(1,0)*{}};
\endxy}

\newcommand{\Hom}{{\rm Hom}}

\newcommand{\Tr}{{\rm Tr}}
\renewcommand{\to}{\rightarrow}

\newcommand{\SU}{{\rm SU}}
\newcommand{\PSU}{{\rm PSU}}
\newcommand{\Br}{{\rm Br}}
\newcommand{\PSL}{{\rm PSL}}

\newcommand{\U}{{\rm U}}
\newcommand{\Lie}{{\rm Lie}}

\def\Id{\mathrm{Id}}

\def\mf{\mathfrak}

\def\Br{{\mathrm{Br}}}

\numberwithin{equation}{section}

\newcommand{\wb}{\overline}

\newcommand{\slt}{{\mathfrak{sl}(2)}}
\newcommand{\Uq}{{U_q\slt}}
\newcommand{\UqMed}{{\wb U_q\slt}}

\newcommand{\UsltH}{{U_q^{H}\slt}}
\newcommand{\Ubar}{{\wb U_q^{H}\slt}}

\let\hat=\widehat

\let\epsilon=\varepsilon

\usepackage{bbm}
\def\C{{\mathbb{C}}}
\def\N{{\mathbbm N}}
\def\R{{\mathbbm R}}
\def\Z{{\mathbbm Z}}

\def\1{\mathbbm{1}}%
\def\ev{\mathrm{ev}}%
\def\coev{\mathrm{coev}}%
\usepackage{bbm}

\newcommand\nc{\newcommand}
\nc\rnc{\renewcommand}
\nc\Kar{\operatorname{Kar}}
\nc\End{\operatorname{End}}

\newcommand{\scs}{\scriptstyle}

\nc\Sym{\operatorname{Sym}}

\allowdisplaybreaks

\newcommand{\tev}{\stackrel{\longleftarrow}{\operatorname{ev}}}
\newcommand{\tcoev}{\stackrel{\longleftarrow}{\operatorname{coev}}}

\newcommand{\qr}{{q}}

\newcommand{\unit}{\ensuremath{\mathbb{I}}}

\newcommand{\qn}[1]{{\left\{#1\right\}}}
\newcommand{\qN}[1]{{\left[#1\right]}}

\newcommand{\qd}{{\mathsf d}}
\newcommand{\qdim}{\operatorname{qdim}}

\newcommand{\Span}{\operatorname{Span}}

\newcommand{\Cp}{{\ddot\C}}

\newcommand{\et}{{\quad\text{and}\quad}}

\newcommand{\Log}{\mathrm{Log}}

\title{Density and unitarity of the Burau representation from a non-semisimple TQFT}

\begin{document}

\author[N. Geer]{Nathan Geer}
\address{Mathematics \& Statistics\\
  Utah State University \\
  Logan, Utah 84322, USA}
  \email{nathan.geer@gmail.com}
\author[A.D. Lauda]{Aaron D. Lauda}
\address{Department of Mathematics \& Department of Physics\\
 University of Southern Californa \\
  Los Angeles, California 90089, USA}
  \email{lauda@usc.edu}

%
\author[B. Patureau-Mirand]{Bertrand Patureau-Mirand}
\address{UMR 6205, LMBA, universit\'e de Bretagne-Sud, universit\'e
  europ\'eenne de Bretagne, BP 573, 56017 Vannes, France }
\email{bertrand.patureau@univ-ubs.fr}
\author[J. Sussan]{Joshua Sussan}
\address{Department of Mathematics\\
  CUNY Medgar Evers \\
  Brooklyn, NY 11225, USA}
  \email{jsussan@mec.cuny.edu}
  \address{Mathematics Program\\
 The Graduate Center, CUNY \\
  New York, NY 10016, USA}

\begin{abstract}
We study the density of the Burau representation from the perspective of a non-semisimple TQFT at a fourth root of unity.  This gives a TQFT construction of Squier's Hermitian form on the Burau representation with possibly mixed signature.  We prove that the image of the braid group in the space of possibly indefinite unitary representations is dense.  We also argue for the potential applications of non-semisimple TQFTs toward topological quantum computation.
  \end{abstract}

\maketitle
\setcounter{tocdepth}{3}

\section{Introduction}

There is a rich interplay between Topological Quantum Field Theory (TQFT) and fault-tolerant approaches to quantum computation via topological quantum computation (TQC)~\cite{NKK,FKLW,FLWmodular}.  In such theories, the coherence of quantum mechanical states is encoded in the topologies of the systems, and unitary operations are performed by braiding of quasiparticles.   Mathematically, these theories are described within the framework of modular tensor categories and have close connections with the representation theory of quantum groups.  

\subsection{Chern-Simons-Witten Theory and TQC}
One of the most studied theories coming from Chern-Simons-Witten theory and its mathematical incarnation via Witten-Reshitkhin-Turaev TQFTs is  associated with
the small quantum group for $\mathfrak{sl}_2$ where the quantum parameter is specialized to a root of unity.  In this framework, quasiparticles correspond to irreducible representations, multiple quasiparticles correspond to a tensor product of irreducibles, and the Hilbert space of the system is formulated from the fusion channels in the tensor product.  This can be identified with various hom spaces within a certain semisimplification of the category of representations of small quantum $\mathfrak{sl}_2$ at this root of unity.

The unitarity of the TQFT then equips this vector space with a Hermitian inner product and the mapping class group of this (2+1)d TQFT induces unitary transformations corresponding to the braiding of quasiparticles.    In the context of topological quantum computation, the critical issue of universality becomes the mathematical question of whether the braid group representations are dense in the corresponding projective unitary group.  There is a great deal of literature studying such questions~\cite{FLWmodular,FLW02}.

\subsection{Non-semisimple TQFT and TQC}
In recent years, there has been a surprising discovery of non-semisimple TQFTs~\cite{BCGP1,BCGP2,CGP1,GPT2}.   These theories are built from representations of quantum trace zero that would have been thrown out in the semisimplification process used in the standard approaches to TQFTs via modular tensor categories.  Already at a fourth root of unity, non-semisimple TQFTs coming from unrolled quantum groups contain information that their semisimple counterparts do not.  For example, it was shown in \cite{BCGP2} that a certain choice of parameters leads to a state space of the torus, which is 2-dimensional, and that the mapping class group action on it is faithful.  This is in stark contrast to the finite image of this mapping class group in the traditional semisimple TQFT.  

Given the critical role of the mapping class group in braiding quasiparticles in quantum computation, we aim to demonstrate how the more sophisticated encoding of topological information in non-semisimple TQFTs can be leveraged to advance TQC.  
This motivated the present work to study the density and efficiency of state spaces for non-semisimple TQFTs already in the simplest case of a fourth root of unity.

%

In our previous work~\cite{GLPMS2,GLPMS3}, we have advocated that the topological advantages of non-semisimple TQFTs may translate into potential advantages for constructing quantum systems based on topological phases of matter. 
In \cite{GLPMS}, we studied non-semisimple TQFTs associated with unrolled $\mf{sl}_2$ and defined a Hermitian structure equipping the tensor spaces with Hermitian forms of possibly mixed signature.  This was extended to more general (super) Lie algebras in ~\cite{GLPMS3}.  

Motivated by the application of non-semisimple TQFTs to problems in topological quantum computation, it is natural to investigate 
\begin{enumerate}
  \item situations in which the Hermitian form arising from non-semisimple TQFTs is positive definite, rather than mixed signature; and
  \item when the associated braid group representations are dense in the associated unitary groups.
\end{enumerate}
As a demonstration that both of these objectives can be achieved in the non-semisimple theory, we investigate the first nontrivial setting by examining the non-semisimple theory of abelian anyons corresponding to taking the quantum parameter at a fourth root of unity\footnote{The semisimple analog would be $\SU(2)_1$ that only supports abelian anyons}. 
  While this case is not interesting from the perspective of quantum computation, it does have applications to classical structures in low-dimensional topology and the study of braid group representations going back to an open problem of Joan Birman from 1974.

\subsection{A quantum setting for the Burau representation }
The theory of $\mathfrak{sl}_2$ at a fourth root of unity is closely connected to the classical Burau representation.   This is a representation of the braid group depending on a parameter $s$.  This well-known representation is deeply connected to numerous areas of mathematics and physics.  Joan Birman posed the question of characterizing the image of the braid group under the (reduced) Burau representation 
\[
\rho \maps B_n \to {\rm GL}_{n-1}(R)
\]
in the general linear group over $R=\Z[s,s^{-1}]$.  This is question 14 in her seminal book ~\cite{Birman}.

There has been some progress towards Birman's question. 
Stoimenow \cite{Stoimenow2006} proved that for certain values of $s$, the closure of the image is dense in the unitary group.  More recently, Nick Salter strongly approximated the answer~ \cite{Salter}.
Salter's solution relied on earlier work of Squier~\cite{Squier}, who defined a Hermitian bilinear pairing on the Burau representation, making the Burau representation unitary.  For a nice exposition of this form and a study of discrete real specializations of the Burau representation, see the work of Scherich \cite{Scherich1, Scherich2}.  Salter approaches Birman's question by proving that the image of the braid group under the Burau representation is $s$-adically dense in the corresponding unitary group; here density is measured in the \emph{$s$-adic topology} where $s$ is the standard parameter in the Burau representation.  In this topology, two matrices $M$ and $N$ are close if $MN^{-1}=\Id$ up to some large power of $s-1$. 

The density of braid group representations has been studied by numerous authors.  Freedman, Larsen, and Wang proved density for these representations when the parameter of the Burau representation is a root of unity~\cite{FLW02}. 
McMullen \cite{McM} studied unitary representations of the braid group coming from the homology of certain Riemann surfaces, and, as a consequence, obtained a bilinear pairing on the Burau representation. He determined when the images of his braid group representations are discrete subgroups.
Venkataramana studied questions about when these subgroups are arithmetic \cite{Venk}.  

\subsection{Unitarity and density from non-semisimple TQFTs}

Here, we show that our prior work ~\cite{GLPMS} studying Hermitian structures in the context of non-semisimple TQFTs gives a completely new construction of the bilinear pairing on the Burau representation from a purely TQFT perspective.  We also use this perspective to prove density results for the images of the representations into both compact and noncompact Lie groups, depending on the parameters.  This includes a strengthening of Salter's result answering Birman's question, where we establish the density of the Burau representation in the appropriate unitary group using the standard topology. 

To state our results more precisely, we study representations of unrolled $\mathfrak{sl}_2$ at a fourth root of unity.   For connections to the Burau representations, we analyze the so-called generic part of the category.  These representations were first studied by Martel \cite{Martel}, who proved faithfulness for the four punctured sphere.
We have a family of representations $V_{\alpha}$ depending upon an irrational real number $\alpha$.  At a fourth root of unity, the tensor product $V_{\alpha} \otimes V_{\alpha}$ decomposes into two irreducibles $V_{2\alpha+1} \oplus V_{2\alpha-1}$. The vector spaces $\mathcal{H}_{n,k,\alpha}$  admitting an action of the braid group $\Br_n$ are formed from $V_{\alpha}^{\otimes n}$ fusing into $V_{n\alpha+k}$.  Our previous work \cite{GLPMS} equips these spaces with Hermitian forms of potentially mixed signature.  We thus obtain braid group actions living in indefinite unitary matrices.
When the form has both positive and mixed signatures, the subgroup obtained is infinite.  This is in contrast to the classical semisimple situation at a fourth root of unity.

One of our main results is the following:
\begin{theo}  
 There is a family of homomorphisms of the braid group $\Br_n$ into $\PSU(\mathcal{H}_{n,k,\alpha})$.
    When $k=\pm(n-3)$, this is the Burau representation (and its dual), where $s=iq^{\alpha}$ is the standard parameter in the Burau representation (Proposition \ref{prop:burauiso}).
    When $\alpha$ is irrational, the image is dense (Theorem \ref{mainthm}).
\end{theo}
We highlight that this proof of density differs from standard strategies to proving of density in the quantum computation literature, see, for example~\cite {FLW02}, as well as Salter's proof of density in the s-adic topology.  Indeed, our proof makes use of a special family of braids known as Jucys-Murphy braids that act diagonally on a natural basis of $\mathcal{H}_{n,k,\alpha}$ and generate a dense subgroup of the maximal torus subgroup of diagonal matrices in the corresponding unitary group. Such a density result is again in stark contrast to the semisimple setting, where such braids would only generate a finite subgroup of diagonal matrices.  

Proposition~\ref{prop:exterior} relates other values of $k$ with exterior powers of the Burau representation.  

\subsection{Efficiency }

Beyond density in the corresponding unitary group, for algorithmic implementation, it is often helpful to study the efficiency of a dense filling.  A foundational result in quantum computation is the Kitaev-Solovay Algorithm \cite{MR1611329}.  This theorem provides an algorithm for approximating a given element in a unitary group from a finite set of elements that densely fills the group.
The number of generators needed grows polylogarithmically in $\frac{1}{\epsilon}$ where $\epsilon$ is the desired precision.  For measure theoretic reasons, the exponent of this polynomial can never be lower than one.  In a remarkable recent work of Kuperberg \cite{Kup-cubic}, the exponent was lowered to $\Log_{\phi} 2$ where $\phi$ is the golden ratio.

We conjecture that the non-semisimple model associated with the generic part of the category of unrolled $\mf{sl}_2$ at a fourth root of unity achieves an optimal level of efficiency for a single qubit (although we do not have an algorithm to produce the approximation).  When $n=3$ and $k=0$, the space  $\mathcal{H}_{3,0,\alpha}$ is two-dimensional and admits a positive definite form for certain values of $\alpha$, see Lemma~\ref{lem:goodalpha}. In Section \ref{sec:quantum}, we numerically show that for many values of $\alpha$, our braiding matrices more efficiently approximate certain gates than the analogous matrices coming from the Fibonacci model.



\subsection{Infinite braid representations from singular subcategory}
We consider one additional model coming from the singular part of the category of representations of unrolled $\mf{sl}_2$ at a fourth root of unity.
Using certain non-semisimple modules, there are natural vector spaces to serve as potential qubit models.  Again, these spaces admit Hermitian forms of possibly mixed signature.  However, we show that the image of the braid group generates a discrete subgroup and thus is not dense in $\PSU(\mathcal{H})$ for some Hermitian vector space $\mathcal{H}$.

\begin{theo} 
 There is an action of the braid group on a Hermitian vector space with mixed signatures coming from the singular part of the category of the unrolled quantum group.  The image into the corresponding indefinite unitary group is infinite but not dense in general (Theorem \ref{thm:sing}).
\end{theo}

\subsection{Organization}
In Section~\ref{sec:burau}, we provide an exposition of the Burau representation of the braid group.
In Section~\ref{sec:unrolled}, we review the unrolled quantum group for $\mathfrak{sl}_2$ at a fourth root of unity and its category of representations.
In Section~\ref{sec:gen}, we study the generic part of the category and prove the key density result. 
Section~\ref{sec:quantum} analyzes the efficiency of the density result from the previous section and provides potential links to quantum computation.
In Section~\ref{sec:sing}, the singular part of the category is studied.  As opposed to the generic part of the category, the corresponding representations of the braid group do not have dense images in indefinite unitary groups and are thus not suitable for quantum computation.
Finally, in Section \ref{app:alt}, we give some alternate approaches to studying density of the braid group representations, which may be of independent interest.

\subsection{Acknowledgements}
The authors would like to thank Sung Kim, Nick Salter, and Paolo Zanardi for helpful remarks.  We are especially grateful to Emmanuel Wagner for pointing out the connection to the Burau representation.  The authors would like to thank Vir Bulchandani and Zachary Stier for very valuable comments on an earlier version of the paper.

N.G.\ is partially
supported by NSF grant DMS-2104497.  A.D.L.\ is
partially supported by NSF grants DMS-1902092 and DMS-2200419, the Army Research Office
W911NF-20-1-0075, and the Simons Foundation collaboration grant on New Structures in Low-dimensional topology.  
J.S.\ is partially supported by a Simons Foundation Travel Support Grant and PSC CUNY Enhanced Award 66685-00 54.
Computations associated with this project were conducted utilizing the Center for Advanced Research Computing (CARC) at the University of Southern California.

\section{The Burau representation} \label{sec:burau}
Let $\Br_n$ be the braid group on $n$ strands.  
That is:
$$\Br_n=\left\langle \sigma_1, \ldots, \sigma_{n-1} | \;
 \begin{array}{ll} 
  \sigma_i \sigma_{i+1} \sigma_i =  \sigma_{i+1} \sigma_i \sigma_{i+1}  ,  \; &i=1,\ldots, n-2   \\
\sigma_i\sigma_j = \sigma_j\sigma_i, &|i-j|>2.  
\end{array}
  \right\rangle $$

In this section, we review the Burau representation of $\Br_n$  following the exposition\footnote{We actually take the inverse of his generators.} of Squier~\cite{Squier}.
The Burau representation of $\Br_n$ on the vector space $B_{n,s}$ spanned by $ \{E_1, \ldots, E_{n-1} \}$ is defined by:
\begin{equation}
    \sigma_i(E_j) = 
    \begin{cases}
        E_i + E_{i-1} & \text{ if } j=i-1, \\
        -s^{-2} E_i & \text{ if } j=i, \\
        s^{-2} E_i + E_{i+1} & \text{ if } j=i+1, \\
        E_j & \text{ otherwise},
    \end{cases}
\end{equation}
where $s$ is a complex number.

Consider a new basis of $B_{n,s}$ given by $ \{f_1, \ldots, f_{n-1} \}$ where
\[
f_j=E_1+(1+s^2)E_2 + (1+s^2+s^4)E_3 + \cdots + (1+s^2+\cdots + s^{2j-2}) E_j .
\]
In this basis, 
$ \sigma_j (f_i) = f_i $ if $i\neq j-1,j$.
In the ordered basis $\{f_{j-1},f_j \}$, the action of this generator is given by
\begin{equation} \label{eq:burii+1}
    \sigma_i = 
    X_i D_B X_i^{-1} 
    \quad \quad
    X_i = \begin{pmatrix}
   s^{i+2}-s^i      &  s^{2i+2}-1 \\
  s^{i}-s^{i+2}      &  s^{2i}-s^{2} 
    \end{pmatrix} 
    \quad \quad
    D_B=
    \begin{pmatrix}
        -s^{-2} & 0 \\
      0  & 1
    \end{pmatrix} \ .
\end{equation} 

Squier found a non-degenerate Hermitian sesquilinear pairing on the Burau representation \cite{Squier}.  In the basis $\{E_1, \ldots, E_{n-1} \}$ it is:
\[
(E_i, E_j) =
\begin{cases}
    s+s^{-1} & \text{ if } j=i, \\
    -s^{i-j} & \text{ if } j=i \pm 1, \\
    0 & \text{ otherwise. }
\end{cases}  
\]

The basis $\{{f}_1, \ldots, {f}_{n-1} \} $ is orthogonal with respect to this form with
\begin{equation} \label{eq:formf_i}
    ({f}_i, {f}_j) = \delta_{i,j} \cdot s^{-2i+1} \frac{(1-s^{2i})(1-s^{2i+2})}{(1-s^2)^2} \ .
\end{equation}

We now come to the foundational result of Squier.

\begin{theorem}[\cite{Squier}]
   If $s$ is on the unit circle, then the Burau representation is unitary.  
\end{theorem}

\section{Unrolled $\mathfrak{sl}_2$} \label{sec:unrolled}
In this section, we recall the algebra $\Ubar$ and a category of
modules over this algebra.  Let $\C$ be
the complex numbers and $\Cp=(\C\setminus \Z)\cup 2\Z.$ Let
$q=e^\frac{\pi\sqrt{-1}}{2}$ be a $4^{th}$-root of unity.  We use the
notation $q^x=e^{\frac{\pi\sqrt{-1} x}{2}}$.  For $n\in \N$, we also set
 $$\qn{x}=q^x-q^{-x},\quad\qN{x}=\frac{\qn x}{\qn1},\quad\qn{n}!=\qn{n}\qn{n-1}\cdots\qn{1}\et\qN{n}!=\qN{n}\qN{n-1}\cdots\qN{1}.$$
Note that for $q=e^\frac{\pi\sqrt{-1}}{2}$ we have the identities
\begin{equation} \label{eq:qidentities}
[\alpha+2] = -[\alpha], \qquad [2-\alpha] = [\alpha]. 
\end{equation}

\subsection{The Drinfel'd-Jimbo quantum group}
Let $\Uq$ be the $\C$-algebra given by generators $E, F, K, K^{-1}$
and relations:
\begin{align}\label{E:RelDCUqsl}
  KK^{-1}&=K^{-1}K=1, & KEK^{-1}&=q^2E, & KFK^{-1}&=q^{-2}F, &
  [E,F]&=\frac{K-K^{-1}}{q-q^{-1}}.
\end{align}
The algebra $\Uq$ is a Hopf algebra where the coproduct, count, and
antipode are defined by
\begin{align}\label{E:HopfAlgDCUqsl}
  \Delta(E)&= 1\otimes E + E\otimes K,
  &\varepsilon(E)&= 0,
  &S(E)&=-EK^{-1},
  \\
  \Delta(F)&=K^{-1} \otimes F + F\otimes 1,
  &\varepsilon(F)&=0,& S(F)&=-KF,
    \\
  \Delta(K)&=K\otimes K
  &\varepsilon(K)&=1,
  & S(K)&=K^{-1}
.\label{E:HopfAlgDCUqsle}
\end{align}
Let $\UqMed$ be the algebra $\Uq$ modulo the relations
$E^2=F^2=0$.

 \subsection{A modified version of $\Uq$}\label{SS:UqH}  Let $\UsltH$ be the
$\C$-algebra given by generators $E, F, K, K^{-1}, H$ and
relations in \eqref{E:RelDCUqsl} along with the relations:
\begin{align*}
  HK&=KH,
& [H,E]&=2E, & [H,F]&=-2F.
\end{align*}
The algebra $\UsltH$ is a Hopf algebra where the coproduct, counit, and
antipode are defined in
\eqref{E:HopfAlgDCUqsl}--\eqref{E:HopfAlgDCUqsle} and by
\begin{align*}
  \Delta(H)&=H\otimes 1 + 1 \otimes H,
  & \varepsilon(H)&=0,
  &S(H)&=-H.
\end{align*}
Define \emph{the unrolled quantum group} $\Ubar$ to be the Hopf algebra $\UsltH$ modulo the relations
$E^2=F^2=0$.

Let $V$ be a finite-dimensional $\Ubar$-module.  An eigenvalue
$\lambda\in \C$ of the operator $H:V\to V$ is called a \emph{weight}
of $V$ and the associated eigenspace is called a \emph{weight space}.
A vector $v$ in the $\lambda$-eigenspace
of $H$ is a \emph{weight vector} of \emph{weight} $\lambda$, i.e. $Hv=\lambda
v$.  We call $V$ a \emph{weight module} if $V$ splits as a direct sum
of weight spaces and $\qr^H=K$ as operators on $V$, i.e., $Kv=q^\lambda
v$ for any vector $v$ of weight $\lambda$.  Let $\cat$ be the category
of finite-dimensional weight $\Ubar$-modules.

Since $\Ubar$ is a Hopf algebra, $\cat$ is a tensor category where
the unit $\unit$ is the 1-dimensional trivial module $\C$.  Moreover,
$\cat$ is $\C$-linear: hom-sets are $\C$-modules, the composition and
tensor product of morphisms are $\C$-bilinear, and
$\End_\cat(\unit)=\C\Id_\unit$.  When it is clear, we denote the unit
$\unit$ by $\C$.  We say a module $V$ is \emph{simple} if it has no
proper submodules.
For a module $V$ and a morphism $f\in\End_\cat(V)$, we write
$\brk f_V=\lambda\in\C$ if $f-\lambda\Id_V$ is nilpotent.  If $V$ is
simple, then Schur's lemma implies that $\End_\cat(V)=\C\Id_V$. Thus for
$f\in \End_\cat(V)$, we have $f=\brk{f}_V \Id_V$.


We will now recall the fact that the category $\cat$ is a ribbon category.
  Let $V$ and $W$ be
objects of $\cat$.  Let $\{v_i\}$ be a basis of $V$ and $\{v_i^*\}$ be
a dual basis of $V^*=\Hom_\C(V,\C)$.  Then
\begin{align*}
  \coev_V \maps & \C \rightarrow V\otimes V^{*}, \text{ given by } 1 \mapsto \sum
  v_i\otimes v_i^*,  &
  \ev_V \maps & V^*\otimes V\rightarrow \C, \text{ given by }
  f\otimes w \mapsto f(w)
\end{align*}
are duality morphisms of $\cat$.
In \cite{Oh}, Ohtsuki truncates the usual formula of the $h$-adic
quantum $\slt$ $R$-matrix to define an operator on $V\otimes W$ by
\begin{equation}
  \label{eq:R}
  R=\qr^{H\otimes H/2} \sum_{n=0}^{1} \frac{\{1\}^{2n}}{\{n\}!}\qr^{n(n-1)/2}
  E^n\otimes F^n.
\end{equation}
where $q^{H\otimes H/2}$ is the operator given by
$$q^{H\otimes H/2}(v\otimes v') =q^{\lambda \lambda'/2}v\otimes v'$$
for weight vectors $v$ and $v'$ of weights of $\lambda$ and
$\lambda'$. The $R$-matrix is not an element in $\Ubar\otimes \Ubar$.
However the action of $R$ on the tensor product of two objects of
$\cat$ is a well-defined linear map.
Moreover, $R$ gives rise to a braiding $c_{V,W} \maps V\otimes W
\rightarrow W \otimes V$ on $\cat$ defined by $v\otimes w \mapsto
\tau(R(v\otimes w))$ where $\tau$ is the permutation $x\otimes
y\mapsto y\otimes x$.
This braiding follows from the invertibility of the $R$-matrix.  An explicit inverse
(see \cite[Section 2.1.2]{BDGG} and \cite{Oh}) is given by
\begin{equation}
  \label{eq:Rinverse}
  R^{-1}= (\sum_{n=0}^{1} (-1)^n  \frac{\{1\}^{2n}}{\{n\}!}\qr^{-n(n-1)/2}
  E^n\otimes F^n) \qr^{-H\otimes H/2}.
\end{equation}

Let $\theta$ be the operator given by
\begin{equation}
\theta=K^{}\sum_{n=0}^{1}
\frac{\{1\}^{2n}}{\{n\}!}\qr^{n(n-1)/2} S(F^n)\qr^{-H^2/2}E^n
\end{equation}
where $q^{-H^2/2}$ is an operator defined on a weight vector $v_\lambda$ by
$q^{-H^2/2}.v_\lambda = q^{-\lambda^2/2}v_\lambda.$
Ohtsuki shows that the family of maps $\theta_V:V\rightarrow V$ in
$\cat$ defined by $v\mapsto \theta^{-1}v$ is a twist (see
\cite{jM,Oh}).

 Now the ribbon structure on $\cat$ yields right duality morphisms
\begin{equation}\label{E:d'b'}
  \tev_{V}=\ev_{V}c_{V,V^*}(\theta_V\otimes\Id_{V^*})\text{ and }\tcoev_V =(\Id_{V^*}\otimes\theta_V)c_{V,V^*}\coev_V
\end{equation}
which are compatible with the left duality morphisms $\{\coev_V\}_V$ and
$\{\ev_V\}_V$.  These duality morphisms are given explicitly by \begin{align*}
  \tcoev_{V} \maps & \C \rightarrow V^*\otimes V, \text{ where } 1 \mapsto
  \sum v_i^* \otimes K^{}v_i, \\ \tev_{V} \maps & V\otimes V^*\rightarrow
  \C, \text{ where } v\otimes f \mapsto f(K^{-1}v).
\end{align*}
The \emph{quantum dimension} $\qdim(V)$ of an object $V$ in $\cat$ is defined by
\[
\qdim(V)= \brk{\tev_V\circ \coev_V}_\unit=\sum  v_i^*(K^{-1}v_i) \ .
\]

 For $g\in\C/2\Z$, define
$\cat_{g}$ as the full subcategory of weight modules whose weights
are all in the class  $g$ (mod $2\Z$).
Then $\cat=\{\cat_g\}_{g\in \C/2\Z}$ is a $\C/2\Z$-graded category (where
$\C/2\Z$ is an additive group). Let $V\in\cat_g$ and $V'\in\cat_{g'}$.
Then the weights of $V\otimes V'$ are congruent to $g+g' \mod 2\Z$,
and so the tensor product is in $\cat_{g+g'}$.  Also, if $g\neq g'$
 then $\Hom_\cat(V, V')=0$ since morphisms in $\cat$  preserve weights.
Finally, if $f\in V^*=\Hom_\C(V,\C)$, then by definition the action of
$H$ on $f$ is given by $(Hf)(v)=f(S(H)v)=-f(Hv)$
and so
$V^{*}\in\cat_{-g}$.
We call the part of the category where $g=0,1$, the {\it singular} part of the category, and refer to the objects in this part of the category as singular objects. 
The category $\cat_g$ is non-semisimple if $g$ is singular, otherwise $g$ is called {\it generic} and $\cat_g$ is semisimple.

We now consider the following class of finite dimensional highest weight modules.
For each $\alpha\in \C$, we let $V_\alpha$ be the $2$-dimensional
highest weight $\Ubar$-module of highest weight $\alpha + 1$.  The
module $V_\alpha$ has a basis $\{v_0,v_{1}\}$ whose action is
given by
\begin{equation}\label{E:BasisV}
H.v_i=(\alpha + 1-2i) v_i,\quad E.v_i= \frac{\qn i\qn{i-\alpha}}{\qn1^2}
v_{i-1} ,\quad F.v_i=v_{i+1}.
\end{equation}
For all $\alpha\in \C$, the quantum dimension of $V_\alpha$ is zero:
$$\qdim(V_\alpha)= \sum_{i=0}^{1} v_i^*(K^{-1}v_i)=
 \sum_{i=0}^{1} q^{-(\alpha + 1-2i)} =
 q^{-\alpha - 1}\frac{1-q^{4}}{1-q^{2}}=0.$$

 For $a\in \Z$, let $\C^H_{2a}$ be the one
dimensional module in $\cat_{\bar 0}$ where both $E$ and $F$ act by zero and $H$ acts by $2a$. 
Every simple module of $\cat$ is isomorphic to exactly one of the
modules in the list:
\begin{itemize}
	\item  $\C^H_{2a}$, for $a\in \Z$,
	\item  $V_\alpha$  for  $\alpha\in(\C\setminus \Z)\cup
	2 \Z$.
\end{itemize}

Let $P_0$ be the projective and
indecomposable module with highest weight $2$, defined in
Proposition 6.2 of \cite{CGP2}.    Moreover, any indecomposable projective weight module has a
highest weight, and such a module $P\in\cat_{\wb0}\cup\cat_{\wb1}$ with
highest weight $2(k+2)-2$ is isomorphic to $P_0\otimes \C^H_{2k}$.

The simple module $\C^H_{0}$ has an indecomposable projective cover $P_0$ of dimension $4$. A detailed description of this module can be found in \cite[Proposition 6.1]{CGP2}.
A summary can be found in Figure \ref{fig:P_i}.  The vectors $w^Y$, for $Y \in \{R,H,S,L\}$ have weights $2,0,0,-2$ respectively (under the action of $H$).
It is easy to check that $\End(P_0) \cong \C[x]/(x^2)$ where
\begin{equation} \label{def:xmap}
x \colon P_0 \rightarrow P_0 \quad \quad \quad w^H \mapsto w^S .
\end{equation}

\begin{figure}
\[
 \hackcenter{
\begin{tikzpicture}[yscale=-1, scale=0.7,  decoration={markings, mark=at position 0.6 with {\arrow{>}};},]
\draw[postaction={decorate}] (.8,.2) to  (.2,.8);
\draw[postaction={decorate}] (-.8,.2) to  (-.2,.8);
\draw[postaction={decorate}] (-.2,-.8) to  (-.8,-.2);
\draw[postaction={decorate}] (.2,-.8) to  (.8,-.2);
\node at (-1,0) {$\scs w^L$};
\node at (1,0) {$\scs w^R$};
\node at (0,1) {$\scs w^S$};
\node at (0,-1) {$\scs w^H$};
\end{tikzpicture}}
\]
\caption{The weight space structure of the module $P_0$.}
  \label{fig:P_i}
\end{figure}

A key ingredient in the construction of so-called non-semisimple TQFTs from the representations of $\Ubar$ is the notion of a modified trace~\cite{GKP1}.  Taking the modified trace of the identity gives a notion of modified dimension, which can be viewed as renormalizing the representations whose usual quantum dimension is zero. 
The modified dimension of $V_{\alpha}$ is given by 
\[
\qd(V_{\alpha}) = \qd(\alpha) := -\frac{2 \sin \left(\frac{\pi \alpha }{2}\right)}{\sin (\pi  \alpha )} = -\frac{1}{[\alpha+1]}
\]
where the last identity holds because $q$ is a fourth root of unity.

If $\alpha, \beta$, and $\alpha+\beta$ are generic, then it is well known (see for example \cite[Theorem 5.2]{CGP2}) that
\begin{equation} \label{ValphaxVbeta}
V_{\alpha} \otimes V_{\beta} \cong V_{\alpha+\beta+1} \oplus V_{\alpha+\beta-1} \ .
\end{equation}

In the next few sections, we will often use a graphical calculus describing this category of representations described in \cite{CGP1}.
All diagrams are read from bottom to top.

Our main focus in the rest of the paper is the study of representations of the braid group $\Br_n$ on various morphism spaces in the category of representations of the unrolled quantum group.

\section{Unitary representations of the braid group: generic part of the category} \label{sec:gen}
In this section, we will study representations coming from the generic part of the category. Let $\alpha$ be irrational so that $V_{\alpha}$ is a projective simple.
For integers $k,n$ with $n\geq 1$ let $\mathcal{H}_{n,k,\alpha}:=\Hom(V_{n \alpha +k}, V_{\alpha}^{\otimes n})$.
Note that this is non-zero if and only if $k=n-1,n-3,\ldots, 1-n$.
%
We choose specific morphisms 
\begin{equation} \label{eq:Habg}
 \hackcenter{
\begin{tikzpicture}[yscale=-1, scale=0.7,  decoration={markings, mark=at position 0.6 with {\arrow{>}};},]
\draw[very thick, postaction={decorate}] (0,0) to (0,1);
\draw[very thick,  postaction={decorate}] (-.5,-1)to [out=90, in=210] (0,0);
\draw[very thick, postaction={decorate}] (.5,-1) to [out=90, in=-30] (0,0);
\node at (.75,-.8) {$\scs \beta$};
\node at (-.75,-.8) {$\scs \alpha$};
\node at (-.5,.6) {$\scs \gamma$};
\end{tikzpicture} }
\quad \in \quad \Hom(V_{\gamma}, V_{\alpha} \otimes V_{\beta})
\end{equation}
coming from the conventions of \cite{CGP1, CM}.
The vector spaces $\mathcal{H}_{n,k,\alpha}$ can be equipped with a non-degenerate Hermitian pairing, which will be utilized throughout this paper.  For details on this construction, see \cite{GLPMS}.

\begin{lemma}
The dimension of $\mathcal{H}_{n,k,\alpha}$ is $\binom{n-1}{\frac{k+n-1}{2}}$ if $k+n-1$ is even, and is zero otherwise.
\end{lemma}

\begin{proof}
This is a straightforward check using the tensor structure in the category.
See for example \eqref{ValphaxVbeta} or \cite[Theorem 5.2]{CGP2}.
\end{proof}

It is useful to have the notion of an $\mathcal{H}_{n,k,\alpha}$ path.  This is a path ${\bf p}=( \alpha=x_0,x_1,\ldots, x_{n-2}, n\alpha+k)$  in Pascal's triangle \eqref{pascal1} from $\alpha$ to $n \alpha+k$ by taking southeast or southwest steps.

\begin{equation} \label{pascal1}
 \hackcenter{
\begin{tikzpicture}[yscale=1, scale=0.7,  decoration={markings, mark=at position 0.6 with {\arrow{>}};},]
\node at (0,0) {$\scs \alpha $};
\node at (-1,-1) {$\scs 2\alpha-1$};
\node at (1,-1) {$\scs 2\alpha+1$};
\node at (-2,-2) {$\scs 3\alpha-2$};
\node at (0,-2) {$\scs 3\alpha$};
\node at (2,-2) {$\scs 3\alpha+2$};
\node at (0,-3) {$\scs \cdots$};
\node at (-3,-4) {$\scs n\alpha-(n-1)$};
\node at (-1,-4) {$\scs \cdots$};
\node at (0,-4) {$\scs n\alpha+k$};
\node at (1,-4) {$\scs \cdots$};
\node at (3,-4) {$\scs n\alpha+(n-1)$};
\end{tikzpicture}}
\end{equation}
We will also label such paths by direction sequences ${\bf s}$ of $L$'s and $R$'s, which indicates at each step whether the path is going to the left or right.  These sequences will have $n-1$ entries.
For example, there is the path $(\alpha,2\alpha-1,3\alpha)$ given in \eqref{pascal2}.  The corresponding direction sequence is $(L,R)$.
\begin{equation} \label{pascal2}
 \hackcenter{
\begin{tikzpicture}[yscale=1, scale=0.7,  decoration={markings, mark=at position 0.6 with {\arrow{>}};},]
\node at (0,0) {$\scs \alpha $};
\node at (-1,-1) {$\scs 2\alpha-1$};
\node at (1,-1) {$\scs 2\alpha+1$};
\node at (-2,-2) {$\scs 3\alpha-2$};
\node at (0,-2) {$\scs 3\alpha$};
\node at (2,-2) {$\scs 3\alpha+2$};
\draw[postaction={decorate}] (-.1,-.1) to  (-.9,-.9);
\draw[postaction={decorate}] (-.9,-1.1) to  (-.1,-1.9);
\end{tikzpicture}}
\end{equation}

\begin{lemma} \label{pathbij}
A basis of the vector space $\mathcal{H}_{n,k,\alpha}$ can be enumerated by paths $\alpha$ to $n \alpha + k$.
This basis is also enumerated by direction sequences with $\frac{n-1+k}{2}$ $R$'s and $\frac{n-1-k}{2}$ $L$'s.
\end{lemma}

\begin{proof}
The bijection maps a path ${\bf p}=(\alpha,x_1,\ldots, x_{n-2}, n\alpha+k)$ to the homomorphism encoded by the tree
\begin{equation} \label{tree1}
 \hackcenter{
\begin{tikzpicture}[yscale=-1, scale=0.7,  decoration={markings, mark=at position 0.6 with {\arrow{>}};},]
\draw[very thick, postaction={decorate}] (0,0) to [out=90, in=220] (.75,1);
\draw[very thick, postaction={decorate}] (1.5,-1) to [out=90, in=-30] (.75,1);
\draw[very thick,  postaction={decorate}] (.8,1) to [out=90, in=220]  (1.35,2);
\draw[very thick,  postaction={decorate}] (-.5,-1)to [out=90, in=210] (0,0);
\draw[very thick, postaction={decorate}] (.5,-1) to [out=90, in=-30] (0,0);
\draw[very thick, postaction={decorate}] (2.5,-1) to [out=90, in=-30] (1.35,2);
\draw[very thick, postaction={decorate}] (1.35,2) to [out=90, in=220] (1.45,2.4);
\node at (1.45,2.5) {$\cdots$};
\draw[very thick, postaction={decorate}] (1.55,2.6) to [out=90, in=220] (1.85,3);
\draw[very thick, postaction={decorate}] (3.5,-1) to [out=90, in=-30] (1.85,3);
\draw[very thick, postaction={decorate}] (1.85,3) to [out=90, in=220] (2.35,4);
\node at (.75,-.8) {$\scs V_{\alpha}$};
\node at (-.75,-.8) {$\scs V_{\alpha}$};
\node at (1.75,-.8) {$\scs V_{\alpha}$};
\node at (.6,2) {$\scs V_{x_2}$};
\node at (-.2,.6) {$\scs V_{x_1}$};
\node at (3,-.8) {$\cdots$};
\node at (2.25,-.8) {$\scs V_{\alpha}$};
\node at (3.75,-.8) {$\scs V_{\alpha}$};
\node at (.8,2.6) {$\scs V_{x_{n-2}}$};
\node at (2,4) {$\scs V_{n \alpha + k}$};
\end{tikzpicture} } \ .
\end{equation}
\end{proof}
We will usually abbreviate the label $V_{\beta}$ of an edge of a tree in \eqref{tree1} simply by $\beta$.  

Lemma~\ref{pathbij} gives a basis for $\mathcal{H}_{n,k,\alpha}$ indexed by fusion diagrams representing homomorphisms $V_{n\alpha+k} \to V_{\alpha}^{\otimes n}$ constructed from the generators in \eqref{eq:Habg}.   Different ways of combining the generators \eqref{eq:Habg} with identity morphisms to construct homomorphism  $V_{n\alpha+k} \to V_{\alpha}^{\otimes n}$ give rise to different bases for $\mathcal{H}_{n,k,\alpha}$. 




\begin{lemma} \label{6jformulalemma}
We have the following change of basis formulas for the space of morphisms:
\begin{equation} \label{6jI}
 \hackcenter{
\begin{tikzpicture}[yscale=-1, scale=0.7,  decoration={markings, mark=at position 0.6 with {\arrow{>}};},]
\draw[very thick, postaction={decorate}] (0,0) to [out=90, in=220] (.75,1);
\draw[very thick, postaction={decorate}] (1.5,-1) to [out=90, in=-30] (.75,1);
\draw[very thick,  postaction={decorate}] (.8,1) to  (.8,2);
\draw[very thick,  postaction={decorate}] (-.5,-1)to [out=90, in=210] (0,0);
\draw[very thick, postaction={decorate}] (.5,-1) to [out=90, in=-30] (0,0);
\node at (.75,-.8) {$\scs a$};
\node at (-.75,-.8) {$\scs b$};
\node at (1.75,-.8) {$\scs a$};
\node at (.2,1.6) {$\scs 2a+b$};
\node at (-.55,.6) {$\scs a+b+1$};
\end{tikzpicture} }
~=~
\md(2a+1) [b-1]
 \hackcenter{
\begin{tikzpicture}[xscale=-1, yscale=-1, scale=0.7,  decoration={markings, mark=at position 0.6 with {\arrow{>}};},]
\draw[very thick, postaction={decorate}] (0,0) to [out=90, in=220] (.75,1);
\draw[very thick, postaction={decorate}] (1.5,-1) to [out=90, in=-30] (.75,1);
\draw[very thick,  postaction={decorate}] (.8,1) to  (.8,2);
\draw[very thick,  postaction={decorate}] (-.5,-1)to [out=90, in=210] (0,0);
\draw[very thick, postaction={decorate}] (.5,-1) to [out=90, in=-30] (0,0);
\node at (.75,-.8) {$\scs a$};
\node at (-.75,-.8) {$\scs a$};
\node at (1.75,-.8) {$\scs b$};
\node at (.2,1.6) {$\scs 2a+b$};
\node at (-.5,.6) {$\scs 2a+1$};
\end{tikzpicture} }
~+~
\md(2a-1) [-a-1]
 \hackcenter{
\begin{tikzpicture}[xscale=-1, yscale=-1, scale=0.7,  decoration={markings, mark=at position 0.6 with {\arrow{>}};},]
\draw[very thick, postaction={decorate}] (0,0) to [out=90, in=220] (.75,1);
\draw[very thick, postaction={decorate}] (1.5,-1) to [out=90, in=-30] (.75,1);
\draw[very thick,  postaction={decorate}] (.8,1) to  (.8,2);
\draw[very thick,  postaction={decorate}] (-.5,-1)to [out=90, in=210] (0,0);
\draw[very thick, postaction={decorate}] (.5,-1) to [out=90, in=-30] (0,0);
\node at (.75,-.8) {$\scs a$};
\node at (-.75,-.8) {$\scs a$};
\node at (1.75,-.8) {$\scs b$};
\node at (.2,1.6) {$\scs 2a+b$};
\node at (-.5,.6) {$\scs 2a-1$};
\end{tikzpicture} } \ ,
\end{equation}

\begin{equation} \label{6jII}
 \hackcenter{
\begin{tikzpicture}[yscale=-1, scale=0.7,  decoration={markings, mark=at position 0.6 with {\arrow{>}};},]
\draw[very thick, postaction={decorate}] (0,0) to [out=90, in=220] (.75,1);
\draw[very thick, postaction={decorate}] (1.5,-1) to [out=90, in=-30] (.75,1);
\draw[very thick,  postaction={decorate}] (.8,1) to  (.8,2);
\draw[very thick,  postaction={decorate}] (-.5,-1)to [out=90, in=210] (0,0);
\draw[very thick, postaction={decorate}] (.5,-1) to [out=90, in=-30] (0,0);
\node at (.75,-.8) {$\scs a$};
\node at (-.75,-.8) {$\scs b$};
\node at (1.75,-.8) {$\scs a$};
\node at (.2,1.6) {$\scs 2a+b$};
\node at (-.55,.6) {$\scs a+b-1$};
\end{tikzpicture} }
~=~
\md(2a+1) [2a+b-1]
 \hackcenter{
\begin{tikzpicture}[xscale=-1, yscale=-1, scale=0.7,  decoration={markings, mark=at position 0.6 with {\arrow{>}};},]
\draw[very thick, postaction={decorate}] (0,0) to [out=90, in=220] (.75,1);
\draw[very thick, postaction={decorate}] (1.5,-1) to [out=90, in=-30] (.75,1);
\draw[very thick,  postaction={decorate}] (.8,1) to  (.8,2);
\draw[very thick,  postaction={decorate}] (-.5,-1)to [out=90, in=210] (0,0);
\draw[very thick, postaction={decorate}] (.5,-1) to [out=90, in=-30] (0,0);
\node at (.75,-.8) {$\scs a$};
\node at (-.75,-.8) {$\scs a$};
\node at (1.75,-.8) {$\scs b$};
\node at (.2,1.6) {$\scs 2a+b$};
\node at (-.5,.6) {$\scs 2a+1$};
\end{tikzpicture} }
~+~
\md(2a-1) [-a-1]
 \hackcenter{
\begin{tikzpicture}[xscale=-1, yscale=-1, scale=0.7,  decoration={markings, mark=at position 0.6 with {\arrow{>}};},]
\draw[very thick, postaction={decorate}] (0,0) to [out=90, in=220] (.75,1);
\draw[very thick, postaction={decorate}] (1.5,-1) to [out=90, in=-30] (.75,1);
\draw[very thick,  postaction={decorate}] (.8,1) to  (.8,2);
\draw[very thick,  postaction={decorate}] (-.5,-1)to [out=90, in=210] (0,0);
\draw[very thick, postaction={decorate}] (.5,-1) to [out=90, in=-30] (0,0);
\node at (.75,-.8) {$\scs a$};
\node at (-.75,-.8) {$\scs a$};
\node at (1.75,-.8) {$\scs b$};
\node at (.2,1.6) {$\scs 2a+b$};
\node at (-.5,.6) {$\scs 2a-1$};
\end{tikzpicture} } \ ,
\end{equation}

\begin{equation} \label{6jIII}
 \hackcenter{
\begin{tikzpicture}[yscale=-1, scale=0.7,  decoration={markings, mark=at position 0.6 with {\arrow{>}};},]
\draw[very thick, postaction={decorate}] (0,0) to [out=90, in=220] (.75,1);
\draw[very thick, postaction={decorate}] (1.5,-1) to [out=90, in=-30] (.75,1);
\draw[very thick,  postaction={decorate}] (.8,1) to  (.8,2);
\draw[very thick,  postaction={decorate}] (-.5,-1)to [out=90, in=210] (0,0);
\draw[very thick, postaction={decorate}] (.5,-1) to [out=90, in=-30] (0,0);
\node at (.75,-.8) {$\scs a$};
\node at (-.75,-.8) {$\scs b$};
\node at (1.75,-.8) {$\scs a$};
\node at (0,1.6) {$\scs 2a+b+2$};
\node at (-.55,.6) {$\scs a+b+1$};
\end{tikzpicture} }
~=~
\md(2a+1) [-a-b-2]
 \hackcenter{
\begin{tikzpicture}[xscale=-1, yscale=-1, scale=0.7,  decoration={markings, mark=at position 0.6 with {\arrow{>}};},]
\draw[very thick, postaction={decorate}] (0,0) to [out=90, in=220] (.75,1);
\draw[very thick, postaction={decorate}] (1.5,-1) to [out=90, in=-30] (.75,1);
\draw[very thick,  postaction={decorate}] (.8,1) to  (.8,2);
\draw[very thick,  postaction={decorate}] (-.5,-1)to [out=90, in=210] (0,0);
\draw[very thick, postaction={decorate}] (.5,-1) to [out=90, in=-30] (0,0);
\node at (.75,-.8) {$\scs a$};
\node at (-.75,-.8) {$\scs a$};
\node at (1.75,-.8) {$\scs b$};
\node at (0,1.6) {$\scs 2a+b+2$};
\node at (-.5,.6) {$\scs 2a+1$};
\end{tikzpicture} } \ ,
\end{equation}

\begin{equation} \label{6jIV}
 \hackcenter{
\begin{tikzpicture}[yscale=-1, scale=0.7,  decoration={markings, mark=at position 0.6 with {\arrow{>}};},]
\draw[very thick, postaction={decorate}] (0,0) to [out=90, in=220] (.75,1);
\draw[very thick, postaction={decorate}] (1.5,-1) to [out=90, in=-30] (.75,1);
\draw[very thick,  postaction={decorate}] (.8,1) to  (.8,2);
\draw[very thick,  postaction={decorate}] (-.5,-1)to [out=90, in=210] (0,0);
\draw[very thick, postaction={decorate}] (.5,-1) to [out=90, in=-30] (0,0);
\node at (.75,-.8) {$\scs a$};
\node at (-.75,-.8) {$\scs b$};
\node at (1.75,-.8) {$\scs a$};
\node at (0,1.6) {$\scs 2a+b-2$};
\node at (-.55,.6) {$\scs a+b-1$};
\end{tikzpicture} }
~=~
\md(2a-1) [2a-2]
 \hackcenter{
\begin{tikzpicture}[xscale=-1, yscale=-1, scale=0.7,  decoration={markings, mark=at position 0.6 with {\arrow{>}};},]
\draw[very thick, postaction={decorate}] (0,0) to [out=90, in=220] (.75,1);
\draw[very thick, postaction={decorate}] (1.5,-1) to [out=90, in=-30] (.75,1);
\draw[very thick,  postaction={decorate}] (.8,1) to  (.8,2);
\draw[very thick,  postaction={decorate}] (-.5,-1)to [out=90, in=210] (0,0);
\draw[very thick, postaction={decorate}] (.5,-1) to [out=90, in=-30] (0,0);
\node at (.75,-.8) {$\scs a$};
\node at (-.75,-.8) {$\scs a$};
\node at (1.75,-.8) {$\scs b$};
\node at (0,1.6) {$\scs 2a+b-2$};
\node at (-.5,.6) {$\scs 2a-1$};
\end{tikzpicture} } \ .
\end{equation}

\end{lemma}

\begin{proof}
This is a straightforward but lengthy calculation using \cite[Equation Nj]{CGP1} or a quicker manipulation of \cite[Proof of Proposition 6.18]{BCGP2} where formulas for $6j$ symbols at a fourth root of unity are provided which uses \cite{CM}.
\end{proof}

The space $\mathcal{H}_{n,k,\alpha}$ has a non-degenerate Hermitian sesquilinear pairing \cite{GLPMS}.

\begin{proposition} \label{prop:pairinggen}   
The non-degenerate Hermitian pairing on $\mathcal{H}_{n,k,\alpha}$ is orthogonal with respect to the basis
${\bf p}=(\alpha,x_1,\ldots, x_{n-2}, x_{n-1}=n\alpha+k)$.
On these basis vectors, the pairing is:
\begin{equation}
\langle {\bf p} | {\bf p} \rangle 
=
\prod_{j=0}^{n-2} \frac{[x_j+\alpha][(x_j+\alpha)^{\frac{x_{j+1}-x_j-\alpha+1}{2}}]}{[(x_j+1)^\frac{x_j+\alpha+1-x_{j+1}}{2}]
[(\alpha+1)^\frac{x_j+\alpha+1-x_{j+1}}{2}]}
\qd(V_{n \alpha+k})\ .
\end{equation}
\end{proposition}

\begin{proof}
This is a lengthy but straightforward calculation using the definitions of \cite{GLPMS} where it is proved that there is an anti-linear involution $\dagger$ on the morphisms of the category of Hermitian modules for the unrolled quantum group.
For two morphisms $f,g$ with the same source and target, the pairing is defined to be the modified trace of $f^{\dagger} g$.

The proposition follows from the straightforward calculations:
\begin{equation}
 \hackcenter{
\begin{tikzpicture}[yscale=-1, scale=0.7,  decoration={markings, mark=at position 0.6 with {\arrow{>}};},]
\draw[very thick, postaction={decorate}] (0,0) to (0,1);
\draw[very thick,  postaction={decorate}] (-.5,-1)to [out=90, in=210] (0,0);
\draw[very thick, postaction={decorate}] (.5,-1) to [out=90, in=-30] (0,0);
\node at (.75,-.8) {$\scs \beta$};
\node at (-.75,-.8) {$\scs \alpha$};
\node at (-.9,.6) {$\scs \alpha+\beta+1$};
\end{tikzpicture} }^{\dagger}
=
-[\alpha+\beta+2]
 \hackcenter{
\begin{tikzpicture}[yscale=1, scale=0.7,  decoration={markings, mark=at position 0.6 with {\arrow{<}};},]
\draw[very thick, postaction={decorate}] (0,0) to (0,1);
\draw[very thick,  postaction={decorate}] (-.5,-1)to [out=90, in=210] (0,0);
\draw[very thick, postaction={decorate}] (.5,-1) to [out=90, in=-30] (0,0);
\node at (.75,-.8) {$\scs \beta$};
\node at (-.75,-.8) {$\scs \alpha$};
\node at (-.9,.6) {$\scs \alpha+\beta+1$};
\end{tikzpicture} }
\qquad \quad 
 \hackcenter{
\begin{tikzpicture}[yscale=-1, scale=0.7,  decoration={markings, mark=at position 0.6 with {\arrow{>}};},]
\draw[very thick, postaction={decorate}] (0,0) to (0,1);
\draw[very thick,  postaction={decorate}] (-.5,-1)to [out=90, in=210] (0,0);
\draw[very thick, postaction={decorate}] (.5,-1) to [out=90, in=-30] (0,0);
\node at (.75,-.8) {$\scs \beta$};
\node at (-.75,-.8) {$\scs \alpha$};
\node at (-.9,.6) {$\scs \alpha+\beta-1$};
\end{tikzpicture} }^{\dagger}
=
-\frac{1}{[\alpha+1][\beta+1]}
 \hackcenter{
\begin{tikzpicture}[yscale=1, scale=0.7,  decoration={markings, mark=at position 0.6 with {\arrow{<}};},]
\draw[very thick, postaction={decorate}] (0,0) to (0,1);
\draw[very thick,  postaction={decorate}] (-.5,-1)to [out=90, in=210] (0,0);
\draw[very thick, postaction={decorate}] (.5,-1) to [out=90, in=-30] (0,0);
\node at (.75,-.8) {$\scs \beta$};
\node at (-.75,-.8) {$\scs \alpha$};
\node at (-.9,.6) {$\scs \alpha+\beta-1$};
\end{tikzpicture} } \ .
\end{equation}
\end{proof}

\begin{lemma} \label{Br2lemma} \cite[Equation Ng]{CGP1}
There is an action of the braid group $\Br_2$ on the 1-dimensional morphism spaces $\Hom(V_{\alpha+\beta-1},V_{\alpha} \otimes V_{\beta})$ and $\Hom(V_{\alpha+\beta+1},V_{\alpha} \otimes V_{\beta})$ given by:
\[ 
 \hackcenter{
\begin{tikzpicture}[yscale=-1, scale=0.7,  decoration={markings, mark=at position 0.6 with {\arrow{>}};},]
\draw[very thick, postaction={decorate}] (0,0) to (0,1);
\draw[very thick,  postaction={decorate}] (-.5,-1)to [out=90, in=210] (0,0);
\draw[very thick, postaction={decorate}] (.5,-1) to [out=90, in=-30] (0,0);
\node at (.75,-.8) {$\scs \beta$};
\node at (-.75,-.8) {$\scs \alpha$};
\node at (-.9,.6) {$\scs \alpha+\beta-1$};
       \draw [black, very thick]    (-.5,-2) .. controls +(0,.25) and +(0,-.25) ..  (.5,-1) ;
      \path [fill=white] (-.25,-1.4) rectangle (.75,-1.6);
       \draw [black, very thick]   (.5,-2)  .. controls +(0,.25) and +(0,-.25) ..  (-.5,-1);
\end{tikzpicture} }
~=~
q^{\frac{-\alpha-\beta+\alpha \beta+1}{2}}
 \hackcenter{
\begin{tikzpicture}[yscale=-1, scale=0.7,  decoration={markings, mark=at position 0.6 with {\arrow{>}};},]
\draw[very thick, postaction={decorate}] (0,0) to (0,1);
\draw[very thick,  postaction={decorate}] (-.5,-1)to [out=90, in=210] (0,0);
\draw[very thick, postaction={decorate}] (.5,-1) to [out=90, in=-30] (0,0);
\node at (.75,-.8) {$\scs \alpha$};
\node at (-.75,-.8) {$\scs \beta$};
\node at (-.9,.6) {$\scs \alpha+\beta-1$};
\end{tikzpicture} }
 \ ,
\qquad 
 \hackcenter{
\begin{tikzpicture}[yscale=-1, scale=0.7,  decoration={markings, mark=at position 0.6 with {\arrow{>}};},]
\draw[very thick, postaction={decorate}] (0,0) to (0,1);
\draw[very thick,  postaction={decorate}] (-.5,-1)to [out=90, in=210] (0,0);
\draw[very thick, postaction={decorate}] (.5,-1) to [out=90, in=-30] (0,0);
\node at (.75,-.8) {$\scs \beta$};
\node at (-.75,-.8) {$\scs \alpha$};
\node at (-.9,.6) {$\scs \alpha+\beta+1$};
       \draw [black, very thick]    (-.5,-2) .. controls +(0,.25) and +(0,-.25) ..  (.5,-1) ;
      \path [fill=white] (-.25,-1.4) rectangle (.75,-1.6);
       \draw [black, very thick]   (.5,-2)  .. controls +(0,.25) and +(0,-.25) ..  (-.5,-1);
\end{tikzpicture} }
~=~
q^{\frac{\alpha+\beta+\alpha \beta+1}{2}}
 \hackcenter{
\begin{tikzpicture}[yscale=-1, scale=0.7,  decoration={markings, mark=at position 0.6 with {\arrow{>}};},]
\draw[very thick, postaction={decorate}] (0,0) to (0,1);
\draw[very thick,  postaction={decorate}] (-.5,-1)to [out=90, in=210] (0,0);
\draw[very thick, postaction={decorate}] (.5,-1) to [out=90, in=-30] (0,0);
\node at (.75,-.8) {$\scs \alpha$};
\node at (-.75,-.8) {$\scs \beta$};
\node at (-.9,.6) {$\scs \alpha+\beta+1$};
\end{tikzpicture} }
 \ .
\]
\end{lemma}

These formulas could be written uniformly as
\begin{equation} \label{eq:twistuniform}
 \hackcenter{
\begin{tikzpicture}[yscale=-1, scale=0.7,  decoration={markings, mark=at position 0.6 with {\arrow{>}};},]
\draw[very thick, postaction={decorate}] (0,0) to (0,1);
\draw[very thick,  postaction={decorate}] (-.5,-1)to [out=90, in=210] (0,0);
\draw[very thick, postaction={decorate}] (.5,-1) to [out=90, in=-30] (0,0);
\node at (.75,-.8) {$\scs \beta$};
\node at (-.75,-.8) {$\scs \alpha$};
\node at (-.4,.6) {$\scs \gamma$};
       \draw [black, very thick]    (-.5,-2) .. controls +(0,.25) and +(0,-.25) ..  (.5,-1) ;
      \path [fill=white] (-.25,-1.4) rectangle (.75,-1.6);
       \draw [black, very thick]   (.5,-2)  .. controls +(0,.25) and +(0,-.25) ..  (-.5,-1);
\end{tikzpicture} }
~=~
q^{\frac{\gamma^2-\alpha^2-\beta^2+1}{4}}
 \hackcenter{
\begin{tikzpicture}[yscale=-1, scale=0.7,  decoration={markings, mark=at position 0.6 with {\arrow{>}};},]
\draw[very thick, postaction={decorate}] (0,0) to (0,1);
\draw[very thick,  postaction={decorate}] (-.5,-1)to [out=90, in=210] (0,0);
\draw[very thick, postaction={decorate}] (.5,-1) to [out=90, in=-30] (0,0);
\node at (.75,-.8) {$\scs \alpha$};
\node at (-.75,-.8) {$\scs \beta$};
\node at (-.4,.6) {$\scs \gamma$};
\end{tikzpicture} }
 \ .
\end{equation}

We continue with the following easy lemma.
\begin{lemma} \label{bridge1}
Consider paths
\begin{align*}
{\bf p} &=(\alpha,x_1,\ldots,x_{i-1},x_i,x_{i+1},\ldots, x_{n-2}, n\alpha+k) \\
{\bf p}' &=(\alpha,x_1,\ldots,x_{i-1},x_i',x_{i+1},\ldots, x_{n-2}, n\alpha+k).
\end{align*}
If the morphisms 
\begin{equation}\label{E:Tree-xi}
 \hackcenter{
\begin{tikzpicture}[yscale=-1, scale=0.7,  decoration={markings, mark=at position 0.6 with {\arrow{>}};},]
\draw[very thick, postaction={decorate}] (0,0) to [out=90, in=220] (.75,1);
\draw[very thick, postaction={decorate}] (1.5,-1) to [out=90, in=-30] (.75,1);
\draw[very thick,  postaction={decorate}] (.8,1) to  (.8,2);
\draw[very thick,  postaction={decorate}] (-.5,-1)to [out=90, in=210] (0,0);
\draw[very thick, postaction={decorate}] (.5,-1) to [out=90, in=-30] (0,0);
\node at (.75,-.8) {$\scs \alpha$};
\node at (-.95,-.8) {$\scs x_{i-1}$};
\node at (1.75,-.8) {$\scs \alpha$};
\node at (.3,1.6) {$\scs x_{i+1}$};
\node at (-.5,.6) {$\scs x_i$};
\end{tikzpicture} } \ ,
\quad \quad
 \hackcenter{
\begin{tikzpicture}[yscale=-1, scale=0.7,  decoration={markings, mark=at position 0.6 with {\arrow{>}};},]
\draw[very thick, postaction={decorate}] (0,0) to [out=90, in=220] (.75,1);
\draw[very thick, postaction={decorate}] (1.5,-1) to [out=90, in=-30] (.75,1);
\draw[very thick,  postaction={decorate}] (.8,1) to  (.8,2);
\draw[very thick,  postaction={decorate}] (-.5,-1)to [out=90, in=210] (0,0);
\draw[very thick, postaction={decorate}] (.5,-1) to [out=90, in=-30] (0,0);
\node at (.75,-.8) {$\scs \alpha$};
\node at (-.95,-.8) {$\scs x_{i-1}$};
\node at (1.75,-.8) {$\scs \alpha$};
\node at (.3,1.6) {$\scs x_{i+1}$};
\node at (-.5,.6) {$\scs x_i'$};
\end{tikzpicture} } 
\end{equation}
span $\Hom(V_{x_{i+1}},V_{x_{i-1}}\otimes V_{\alpha}\otimes V_{\alpha})$, then the braiding operator $\sigma_{i+1}$ preserves the subspace spanned by ${\bf p}$ and ${\bf p}'$.
\end{lemma}
\begin{proof}
By Lemma \ref{pathbij}, the paths ${\bf p}$ and ${\bf p}'$   correspond to two trees which are identical except in the portions indicated in Equation \eqref{E:Tree-xi}.  

The operator $\sigma_{i+1}$ braids the second and third strands of these trees.  The action on the spanning set \eqref{E:Tree-xi} can be calculated explicitly using Lemmas~\ref{6jformulalemma} and \ref{Br2lemma}, and 
it easily follows that $\sigma_{i+1}$ preserves the subspace spanned by the paths ${\bf p}$ and ${\bf p}'$.
\end{proof}

We now restrict ourselves to the special case $\mathcal{H}_{n,n-3,\alpha}$.
   Consider the basis of direction sequences
    \begin{equation} \label{eq:eibasis}
        e_i=(R, \ldots, R, L, R, \ldots, R)
    \end{equation}
    where the $L$ is in position $i$. 
    One easily calculates that that in this basis $\sigma_1(e_j)=q^{\frac{(\alpha+1)^2}{2}} e_j$,  except $\sigma_1(e_1)=q^{\frac{(\alpha-1)^2}{2}} e_1$.
    Furthermore, $\sigma_i(e_j)=q^{\frac{(\alpha+1)^2}{2}} e_j$ if $i\geq 2$ and $j \neq i-1,i$. Finally, in the ordered basis $\{e_{i-1},e_i\}$, the action of this generator is given by 
    \begin{equation} \label{eq:actionn-3}
      \sigma_i = A_i^{-1} \sigma A_i 
    \end{equation}
    where 
    \begin{equation} \label{eq:defofsigma}
        \sigma= \begin{pmatrix}
           q^{\frac{(\alpha-1)^2}{2}} & 0 \\
           0 & q^{\frac{(\alpha+1)^2}{2}}
        \end{pmatrix} \ ,
        \quad 
        A_i = 
        \begin{pmatrix}
            \md(2\alpha-1) [-\alpha-1] & \md(2\alpha-1) [-\alpha-1] \\
            \md(2\alpha+1)[(i+1) \alpha + i-3] &
            \md(2\alpha+1)[(i-1) \alpha + i-3]
        \end{pmatrix} \ .
    \end{equation}
In the basis $\{e_1,\ldots,e_{n-1} \}$, the Jucys-Murphy elements (see \eqref{eq:JMdef} later on for a graphical definition of these elements) are given by
\begin{equation} \label{eq:JMn-3}
J_j =
    {\rm Diag} 
    \hspace{.05in}
    q^{j(\alpha+1)^2}\left\{\underbrace{q^{(-2\alpha-2)}, q^{(-2\alpha-2)}, \ldots, 
q^{(-2\alpha-2)}}_{j-1}, q^{(-2j-2)\alpha-(2j-2)}, \underbrace{1, \ldots, 1}_{n-1-j}
\right\} \ .
\end{equation}
\begin{proposition} \label{prop:burauiso}
%
There is an isomorphism of projective representations of the braid group $\Br_n$ from the Burau representation  $B_{n,s}$ from Section~\ref{sec:burau} to the morphism space $\mathcal{H}_{n,n-3,\alpha}$, where $s=iq^{\alpha}$.
Furthermore, up to an overall sign, the isomorphism of projective representations respects the bilinear forms up to a sign determined by the sign of $[\alpha+1][n\alpha+n-2]$.
\end{proposition}

\begin{proof}
  Define an isomorphism of vector spaces $\phi \colon B_{n,s} \rightarrow \mathcal{H}_{n,n-3,\alpha} $  where 
    \begin{equation} \label{eq:burauiso}
  \phi(f_k)=(-1)^{k+1}
  \frac
  {((iq^{\alpha})^{k+1}-(iq^{\alpha})^{-k-1})
  ((iq^{\alpha})^{k}-(iq^{\alpha})^{-k})}
  {((iq^{\alpha})^{2}-(iq^{\alpha})^{-2})
  ((iq^{\alpha})^{}-(iq^{\alpha})^{-1})}
  \left(
  \frac
{[2\alpha+2] \sqrt{[\alpha+1][n\alpha+n-2]}}  
{[2\alpha][3\alpha+1]\cdots[n\alpha+n-2]} 
\right)
  e_k
  \end{equation}
  for $k=1,\ldots,n-1$.
  The second factor in the parenthesis is just needed to normalize the isomorphism of representations so that it becomes an isomorphism of unitary representations.

  Consider the matrix $A_j$ in \eqref{eq:defofsigma}.  Up to a scalar, it is equal to
  \begin{equation}
      \hat{A}_j =
      \begin{pmatrix}
          [\alpha+1] & [\alpha+1] \\
          [(j+1)\alpha+j-3] & [(j-1)\alpha+j-3]
      \end{pmatrix} \ .
  \end{equation}
  Then one calculates
  \begin{equation}
      (\hat{A}_j)^{-1} =
      \begin{pmatrix}
          [(j-1)\alpha+j-3] & -[\alpha+1] \\
          - [(j+1)\alpha+j-3] & [\alpha+1] 
      \end{pmatrix} \ .
  \end{equation}
  Letting $s=iq^{\alpha}$, one obtains that 
   \begin{equation}
      (\hat{A}_j)^{-1} = c
      \begin{pmatrix}
          s^{2}-s^{2j} & s^{j}-s^{j+2} \\
          1-s^{2j+2} & s^{j+2}-s^{j}
      \end{pmatrix} 
  \end{equation}
  for some scalar $c$.
  Thus, up to a scalar, $\sigma_j$ acts on the basis elements
  $\{e_{j-1},e_j \}$
  by
  \begin{equation}
   (\hat{A}_j)^{-1} 
   \begin{pmatrix}
     -s^{-2} & 0 \\
   0  & 1
   \end{pmatrix}
   \hat{A}_j \ .
  \end{equation}
  Next, notice that \eqref{eq:burii+1} could be written as:
  \begin{equation}
 \sigma_j= X_j \Gamma_j D_b \Gamma_j^{-1} X_j^{-1}, \quad \quad 
       \Gamma_j = (-1)^{j+1}
       \begin{pmatrix}
           \frac{s^2-s^{-2}}{s^j-s^{-j}}  & 0 \\
         0   &  \frac{s^2-s^{-2}}{(s^{j-1}-s^{-(j-1)})(s^{j}-s^{-j})(s^{j+1}-s^{-(j+1)})}
       \end{pmatrix} \ .
  \end{equation}
  Clearly the diagonal matrices $D_b$ and $ \begin{pmatrix}
     -s^{-2} & 0 \\
 0    & 1
   \end{pmatrix}$ are the same.
Letting $\phi_{j-1,j}$ be the diagonal matrix containing entries $j-1$ and $j$ of the map $\phi$, it is straightforward to show that
$\phi_{j-1,j} X_j \Gamma_j=(\hat{A}_j)^{-1} $ up to a scalar, so $\phi$ is a projective isomorphism of representations.

  Note that the pairing on $\mathcal{H}_{n,n-3,\alpha}$ coming from Proposition \ref{prop:pairinggen} is
  \[
\langle e_k | e_k \rangle 
= \frac
{[(k+1) \alpha +k-1] \prod_{j=1}^{k-1} [(j+1) \alpha +j-1]^2
\prod_{j=k+1}^{n-1} [(j+1) \alpha + j-3]^2}
{[k\alpha+k][\alpha+1]} \qd(V_{n\alpha+n-3})\ .
  \]
Compare this to Squier's form on the Burau representation from \eqref{eq:formf_i}.
After substituting $s=iq^{\alpha}$, it is straightforward to show that $\phi$ is an isometry.
  \end{proof}

\begin{proposition} \label{prop:exterior}
    There is an isomorphism of projective representations of the braid group 
    \[\Lambda^k \mathcal{H}_{n,n-3,\alpha} \cong \mathcal{H}_{n,n-1-2k,\alpha} .\]
\end{proposition} 

\begin{proof}
Consider the linear isomorphism $\phi \colon \Lambda^k \mathcal{H}_{n,n-3,\alpha} \rightarrow \mathcal{H}_{n,n-1-2k,\alpha} $ mapping
$e_{l_1} \wedge \cdots \wedge e_{l_k}$ to a sequence ${\bf s}$ with $L$'s in positions $l_1,\ldots,l_k$, with $l_1 < \cdots < l_k$.

First consider the braid group action on $\mathcal{H}_{n,n-1-2k,\alpha}$.  For the action of $\sigma_i$, the relevant part of the tree is the part containing the $(i-1)$st and $i$th leaves.  For such a tree, let $\gamma$ be the number of $R$'s minus the number of $L$'s in the corresponding direction sequence in the first $i-2$ entries.

If a direction sequence ${\bf s}$ has $L$'s in entries $i-1$ and $i$, then the relevant part of the corresponding tree is
\[
 \hackcenter{
\begin{tikzpicture}[yscale=-1, scale=0.7,  decoration={markings, mark=at position 0.6 with {\arrow{>}};},]
\draw[very thick, postaction={decorate}] (0,0) to [out=90, in=220] (.75,1);
\draw[very thick, postaction={decorate}] (1.5,-1) to [out=90, in=-30] (.75,1);
\draw[very thick,  postaction={decorate}] (.8,1) to  (.8,2);
\draw[very thick,  postaction={decorate}] (-.5,-1)to [out=90, in=210] (0,0);
\draw[very thick, postaction={decorate}] (.5,-1) to [out=90, in=-30] (0,0);
\node at (.75,-.8) {$\scs \alpha$};
\node at (-1.65,-.8) {$\scs (i-1) \alpha+\gamma$};
\node at (1.75,-.8) {$\scs \alpha$};
\node at (-.6,1.6) {$\scs (i+1) \alpha + \gamma-2$};
\node at (-.7,.6) {$\scs i \alpha+\gamma-1$};
\end{tikzpicture} } \ .
\]
In this case, $\sigma_i$ acts on this basis element by the scalar $q^{\frac{(\alpha-1)^2}{2}}$.

If a direction sequence ${\bf s}$ has $R$'s in entries $i-1$ and $i$, then the relevant part of the corresponding tree is
\[
 \hackcenter{
\begin{tikzpicture}[yscale=-1, scale=0.7,  decoration={markings, mark=at position 0.6 with {\arrow{>}};},]
\draw[very thick, postaction={decorate}] (0,0) to [out=90, in=220] (.75,1);
\draw[very thick, postaction={decorate}] (1.5,-1) to [out=90, in=-30] (.75,1);
\draw[very thick,  postaction={decorate}] (.8,1) to  (.8,2);
\draw[very thick,  postaction={decorate}] (-.5,-1)to [out=90, in=210] (0,0);
\draw[very thick, postaction={decorate}] (.5,-1) to [out=90, in=-30] (0,0);
\node at (.75,-.8) {$\scs \alpha$};
\node at (-1.65,-.8) {$\scs (i-1) \alpha+\gamma$};
\node at (1.75,-.8) {$\scs \alpha$};
\node at (-.6,1.6) {$\scs (i+1) \alpha + \gamma+2$};
\node at (-.7,.6) {$\scs i \alpha+\gamma+1$};
\end{tikzpicture} } \ .
\]
In this case, $\sigma_i$ acts on this basis element by the scalar $q^{\frac{(\alpha+1)^2}{2}}$.

Now consider the case where we have two basis elements whose direction sequences are the same in all entries except the $(i-1)$st and $i$th entries.  In the first of the two basis elements, the sequence has $L$'s and $R$'s in entries $i-1$ and $i$ respectively. The other basis element has $R$ and $L$ in entries $i-1,i$.  The corresponding trees are
\begin{equation} \label{eq:w1w2}
w_1=
 \hackcenter{
\begin{tikzpicture}[yscale=-1, scale=0.7,  decoration={markings, mark=at position 0.6 with {\arrow{>}};},]
\draw[very thick, postaction={decorate}] (0,0) to [out=90, in=220] (.75,1);
\draw[very thick, postaction={decorate}] (1.5,-1) to [out=90, in=-30] (.75,1);
\draw[very thick,  postaction={decorate}] (.8,1) to  (.8,2);
\draw[very thick,  postaction={decorate}] (-.5,-1)to [out=90, in=210] (0,0);
\draw[very thick, postaction={decorate}] (.5,-1) to [out=90, in=-30] (0,0);
\node at (.75,-.8) {$\scs \alpha$};
\node at (-1.65,-.8) {$\scs (i-1) \alpha+\gamma$};
\node at (1.75,-.8) {$\scs \alpha$};
\node at (-.6,1.6) {$\scs (i+1) \alpha + \gamma$};
\node at (-.7,.6) {$\scs i \alpha+\gamma-1$};
\end{tikzpicture} } 
\quad \quad 
w_2=
\hackcenter{
\begin{tikzpicture}[yscale=-1, scale=0.7,  decoration={markings, mark=at position 0.6 with {\arrow{>}};},]
\draw[very thick, postaction={decorate}] (0,0) to [out=90, in=220] (.75,1);
\draw[very thick, postaction={decorate}] (1.5,-1) to [out=90, in=-30] (.75,1);
\draw[very thick,  postaction={decorate}] (.8,1) to  (.8,2);
\draw[very thick,  postaction={decorate}] (-.5,-1)to [out=90, in=210] (0,0);
\draw[very thick, postaction={decorate}] (.5,-1) to [out=90, in=-30] (0,0);
\node at (.75,-.8) {$\scs \alpha$};
\node at (-1.65,-.8) {$\scs (i-1) \alpha+\gamma$};
\node at (1.75,-.8) {$\scs \alpha$};
\node at (-.6,1.6) {$\scs (i+1) \alpha + \gamma$};
\node at (-.7,.6) {$\scs i \alpha+\gamma+1$};
\end{tikzpicture} } \ .
\end{equation}
One calculates using Lemma \ref{6jformulalemma} that $\sigma_i$ acts by the matrix
\[
A_{i,\gamma}^{-1} \sigma A_{i,\gamma}
\]
where
\begin{equation}
A_{i,\gamma}=
\begin{pmatrix}
 \qd(2\alpha-1) [-\alpha-1]
& \qd(2\alpha-1) [-\alpha-1]  \\
    \qd(2\alpha+1) [(i+1)\alpha+\gamma-1]
&   \qd(2\alpha+1) [(i-1)\alpha+\gamma-1]
\end{pmatrix}
\end{equation}
where recall that $\sigma$ is given in \eqref{eq:defofsigma}.

Next we compute the action of $\sigma_i$ on basis elements of $\Lambda^k \mathcal{H}_{n,n-3,\alpha}$.

The basis element in the exterior power corresponding to a direction sequence with $R$'s in entries $i-1,i$ is a wedge
$e_{s_1} \wedge \cdots \wedge e_{s_k}$ where none of the subscripts are $i-1$ or $i$.  In this case using \eqref{eq:actionn-3}, the action of $\sigma_i$ is given by the scalar $q^{\frac{k(\alpha+1)^2}{2}}$.

The basis element in the exterior power corresponding to a direction sequence with $L$'s in entries $i-1,i$ is a wedge
$e_{s_1} \wedge \cdots \wedge e_{s_k}$ where the subscripts $i-1$ and $i$ appear in the expression.
Then $\sigma_i$ acts by the scalar $q^{\frac{(\alpha+1)^2}{2}}$ on $k-2$ of the factors.
Recall that on the basis $\{e_{i-1}, e_i \}$ that $\sigma_i$
acts by the matrix $A_i^{-1} \sigma A_i$. 
A straightforward computation then yields that
\[
\sigma_i(e_{s_1} \wedge \cdots \wedge e_{s_k})
=
q^{\frac{(\alpha-1)^2}{2} + \frac{(k-1)(\alpha+1)^2}{2} }
e_{s_1} \wedge \cdots \wedge e_{s_k} .
\]

Finally, we consider the two basis elements $w_1, w_2$ in \eqref{eq:w1w2}.  These elements correspond to wedges
$e_{s_1} \wedge \cdots \wedge e_{i-1} \wedge e_{s_r} \wedge \cdots \wedge e_{s_k}$ (with $s_r \neq i$) and
$e_{s_1} \wedge \cdots \wedge e_{s_c} \wedge e_i \wedge \cdots \wedge e_{s_k}$ (with $s_c \neq i-1$) respectively.
Then $\sigma_i$ acts by the scalar $q^{\frac{(\alpha+1)^2}{2}}$ on $k-1$ of the factors and by the matrix
$A_i^{-1} \sigma A_i$ on the parts of the wedges 
$e_{i-1} \wedge e_{s_r}$ and $e_{s_c} \wedge e_i$.

Assume for the moment that 
$A_i^{-1} \sigma A_i
=
A_{i,\gamma}^{-1} \sigma A_{i,\gamma}
$.
Then, comparing the actions of $\sigma_i$ on
$ \Lambda^k \mathcal{H}_{n,n-3,\alpha} $ and $\mathcal{H}_{n,n-1-2k,\alpha}$, we see that they are the same up to a factor of $q^{\frac{(k-1)(\alpha+1)^2}{2}}$.
Thus, the actions are the same projectively.

Finally we need to check $A_i^{-1} \sigma A_i
=
A_{i,\gamma}^{-1} \sigma A_{i,\gamma}
$.
The top rows of the matrices $A_i$ and $A_{i,\gamma}$ are the same, so we just need to analyze the bottom row.
The bottom left entries are $[(i+1)\alpha+i-3]$ and 
$[(i+1)\alpha+\gamma-1]$ respectively.
We claim that $\gamma-1$ and $i-3$ differ by an even number.
Let $a_L$ be the number of $L$'s in the first $i-2$ entries and let $a_R$ be the number of $R$'s in the first $i-2$ entries.  Then $a_L+a_R=i-2$ and $-a_L+a_R=\gamma$.  Then, the claim easily follows.
By \eqref{eq:qidentities}, the bottom left entries are the same up to $\pm 1$.
Similarly, the bottom right entries are the same up to $\pm 1$.
Then it follows easily that $A_i^{-1} \sigma A_i
=
A_{i,\gamma}^{-1} \sigma A_{i,\gamma}
$.
\end{proof}

\begin{theorem} \label{mainthm}
The image of the braid group $\Br_n$ in $\PSU(\mathcal{H}_{n,k,\alpha})$, for $k= \pm (n-3)$ is dense, where $\mathcal{H}_{n, k, \alpha} $ is a non-degenerate finite-dimensional Hermitian vector space with possibly a mixed signature.
\end{theorem}

\begin{proof}
Let $\mathcal{H}=\mathcal{H}_{n, k, \alpha} $.
The fact that the representation is (possibly indefinite) unitary follows from \cite[Proposition 5.8]{GLPMS}. 
Throughout this proof we will use basis elements
 ${\bf p}=(x_0=\alpha,x_1,\ldots, x_{n-2}, x_{n-1}=n\alpha+k)$
 from Lemma \ref{pathbij}.

Consider the elements $J_i$, for $i=1,\ldots,n-1$ in the braid group $\Br_n$:
\begin{equation} \label{eq:JMdef}
J_i ~=~
 \hackcenter{
\begin{tikzpicture}[yscale=-1, scale=0.5,  decoration={markings, mark=at position 0.6 with {\arrow{>}};},]
       \draw [black, very thick]    (-.5,-2) .. controls +(0,.25) and +(0,-.25) ..  (.5,-1) ;
      \path [fill=white] (-.25,-1.4) rectangle (.75,-1.6);
       \draw [black, very thick]   (.5,-2)  .. controls +(0,.25) and +(0,-.25) ..  (-.5,-1);
       \draw[very thick, postaction={}] (-1.5,-1) to (-1.5,-2);
       \draw[very thick, postaction={}] (-2.5,-1) to (-2.5,-3);
          \draw [black, very thick]    (-1.5,-3) .. controls +(0,.25) and +(0,-.25) ..  (-.5,-2) ;
      \path [fill=white] (-1.25,-2.4) rectangle (-.25,-2.6);
       \draw [black, very thick]   (-.5,-3)  .. controls +(0,.25) and +(0,-.25) ..  (-1.5,-2);
          \draw[very thick, postaction={}] (.5,-2) to (.5,-6);
          \draw [black, very thick]    (-2.5,-4) .. controls +(0,.25) and +(0,-.25) ..  (-1.5,-3) ;
      \path [fill=white] (-2.25,-3.4) rectangle (-1.25,-3.6);
       \draw [black, very thick]   (-1.5,-4)  .. controls +(0,.25) and +(0,-.25) ..  (-2.5,-3); 
          \draw [black, very thick]    (-2.5,-5) .. controls +(0,.25) and +(0,-.25) ..  (-1.5,-4) ;
      \path [fill=white] (-2.25,-4.4) rectangle (-1.25,-4.6);
       \draw [black, very thick]   (-1.5,-5)  .. controls +(0,.25) and +(0,-.25) ..  (-2.5,-4); 
          \draw[very thick, postaction={}] (-.5,-5) to (-.5,-3);
                    \draw [black, very thick]    (-1.5,-6) .. controls +(0,.25) and +(0,-.25) ..  (-.5,-5) ;
      \path [fill=white] (-1.25,-5.4) rectangle (-.25,-5.6);
       \draw [black, very thick]   (-.5,-6)  .. controls +(0,.25) and +(0,-.25) ..  (-1.5,-5); 
                    \draw [black, very thick]    (-.5,-7) .. controls +(0,.25) and +(0,-.25) ..  (.5,-6) ;
      \path [fill=white] (-.25,-6.4) rectangle (.75,-6.6);
       \draw [black, very thick]   (.5,-7)  .. controls +(0,.25) and +(0,-.25) ..  (-.5,-6); 
          \draw[very thick, postaction={}] (-2.5,-5) to (-2.5,-7); 
          \draw[very thick, postaction={}] (-1.5,-6) to (-1.5,-7); 
          \draw[very thick, postaction={}] (1.5,-1) to (1.5,-7); 
          \node at (1,-4) {$\cdots$};
          \node at (-2,-1.5) {$\cdots$};
                  \node at (-2.5,-.5) {$\scs {1}$}; 
         \node at (-.5,-.5) {$\scs {i}$}; 
            \node at (1.5,-.5) {$\scs {n}$}; 
\end{tikzpicture} } \ .
\end{equation}
These braid elements act diagonally on the basis spanned by trees \eqref{tree1}.  This is clear for $J_1$.  For $J_i$ with $i>1$, this follows
from the fact that strands slide past trivalent vertices as in \eqref{jucyslide}.

\begin{equation} \label{jucyslide}
 \hackcenter{
\begin{tikzpicture}[yscale=-1, scale=0.5,  decoration={markings, mark=at position 0.6 with {\arrow{>}};},]
       \draw [black, very thick]    (-.5,-2) .. controls +(0,.25) and +(0,-.25) ..  (.5,-1) ;
      \path [fill=white] (-.25,-1.4) rectangle (.75,-1.6);
       \draw [black, very thick]   (.5,-2)  .. controls +(0,.25) and +(0,-.25) ..  (-.5,-1);
       \draw[very thick, postaction={}] (-1.5,-1) to (-1.5,-2);
       \draw[very thick, postaction={}] (-2.5,-1) to (-2.5,-3);
          \draw [black, very thick]    (-1.5,-3) .. controls +(0,.25) and +(0,-.25) ..  (-.5,-2) ;
      \path [fill=white] (-1.25,-2.4) rectangle (-.25,-2.6);
       \draw [black, very thick]   (-.5,-3)  .. controls +(0,.25) and +(0,-.25) ..  (-1.5,-2);
          \draw[very thick, postaction={}] (.5,-2) to (.5,-6);
          \draw [black, very thick]    (-2.5,-4) .. controls +(0,.25) and +(0,-.25) ..  (-1.5,-3) ;
      \path [fill=white] (-2.25,-3.4) rectangle (-1.25,-3.6);
       \draw [black, very thick]   (-1.5,-4)  .. controls +(0,.25) and +(0,-.25) ..  (-2.5,-3); 
          \draw [black, very thick]    (-2.5,-5) .. controls +(0,.25) and +(0,-.25) ..  (-1.5,-4) ;
      \path [fill=white] (-2.25,-4.4) rectangle (-1.25,-4.6);
       \draw [black, very thick]   (-1.5,-5)  .. controls +(0,.25) and +(0,-.25) ..  (-2.5,-4); 
          \draw[very thick, postaction={}] (-.5,-5) to (-.5,-3);
                    \draw [black, very thick]    (-1.5,-6) .. controls +(0,.25) and +(0,-.25) ..  (-.5,-5) ;
      \path [fill=white] (-1.25,-5.4) rectangle (-.25,-5.6);
       \draw [black, very thick]   (-.5,-6)  .. controls +(0,.25) and +(0,-.25) ..  (-1.5,-5); 
                    \draw [black, very thick]    (-.5,-7) .. controls +(0,.25) and +(0,-.25) ..  (.5,-6) ;
      \path [fill=white] (-.25,-6.4) rectangle (.75,-6.6);
       \draw [black, very thick]   (.5,-7)  .. controls +(0,.25) and +(0,-.25) ..  (-.5,-6); 
          \draw[very thick, postaction={}] (-2.5,-5) to (-2.5,-7); 
          \draw[very thick, postaction={}] (-1.5,-6) to (-1.5,-7); 
          \draw[very thick, postaction={}] (1.5,-1) to (1.5,-7); 
          \node at (1,-4) {$\cdots$};
          \node at (-2,-1.5) {$\cdots$};
         \node at (-.5,-.5) {$\scs {i}$}; 
          \draw[very thick, postaction={decorate}] (-2,0) to [out=90, in=220] (-1.25,1);
\draw[very thick, postaction={decorate}] (-.5,-1) to [out=90, in=-30] (-1.25,1);
\draw[very thick,  postaction={decorate}] (-1.2,1) to [out=90, in=220]  (-.65,2);
\draw[very thick,  postaction={decorate}] (-2.5,-1)to [out=90, in=210] (-2,0);
\draw[very thick, postaction={decorate}] (-1.5,-1) to [out=90, in=-30] (-2,0);
\draw[very thick, postaction={decorate}] (.5,-1) to [out=90, in=-30] (-.65,2);
\draw[very thick, postaction={decorate}] (-.65,2) to [out=90, in=220] (-.55,2.4);
\node at (-.55,2.5) {$\cdots$};
\draw[very thick, postaction={decorate}] (-.45,2.6) to [out=90, in=220] (-.15,3);
\draw[very thick, postaction={decorate}] (1.5,-1) to [out=90, in=-30] (.15,3);
\draw[very thick, postaction={decorate}] (-.15,3) to [out=90, in=220] (.35,4);
\node at (1,-.8) {$\cdots$};
\end{tikzpicture}} 
\quad = \quad 
 \hackcenter{
\begin{tikzpicture}[yscale=-1, scale=0.5,  decoration={markings, mark=at position 0.6 with {\arrow{>}};},]
\draw[very thick, postaction={}] (0,0) to [out=90, in=220] (.75,1);
\draw[very thick, postaction={decorate}] (2.5,-6) to [out=90, in=-30] (.75,1);
\draw[very thick,  postaction={decorate}] (.8,1) to  (.8,2);
\draw[very thick,  postaction={decorate}] (-.5,-1)to [out=90, in=210] (0,0);
\draw[very thick, postaction={decorate}] (.5,-1) to [out=90, in=-30] (0,0);
\draw [black, very thick]    (-.5,-2) .. controls +(0,.25) and +(0,-.25) ..  (.5,-1) ;
      \path [fill=white] (-.25,-1.4) rectangle (.75,-1.6);
       \draw [black, very thick]   (.5,-2)  .. controls +(0,.25) and +(0,-.25) ..  (-.5,-1); 
          \draw [black, very thick]    (-.5,-3) .. controls +(0,.25) and +(0,-.25) ..  (.5,-2) ;
      \path [fill=white] (-.25,-2.4) rectangle (.75,-2.6);
       \draw [black, very thick]   (.5,-3)  .. controls +(0,.25) and +(0,-.25) ..  (-.5,-2); 
\draw[very thick, postaction={decorate}] (-1.3,-5) to [out=90, in=220] (-.55,-4);
\draw[very thick, postaction={decorate}] (.2,-6) to [out=90, in=-30] (-.55,-4);
\draw[very thick,  postaction={decorate}] (-.5,-4) to  (-.5,-3);
\draw[very thick,  postaction={decorate}] (-1.8,-6)to [out=90, in=210] (-1.3,-5);
\draw[very thick, postaction={decorate}] (-.8,-6) to [out=90, in=-30] (-1.3,-5);
\draw[very thick, postaction={decorate}] (1.5,-6) .. controls +(0,1.75) and +(0,-1) ..  (.5,-3);
      \path [fill=white] (-.25,.4) rectangle (.75,.85);
      \node at (.25,.6) {$\ddots$};
       \node at (-1.25,-6) {$\cdots$};
\end{tikzpicture} } 
\ .
\end{equation}

Let $G$ denote the closure of the image of the representation in the unitary group $\U(\mathcal{H})$, so $G \subset \U(\mathcal{H})$. 
$G$ is a Lie group and let $\mathfrak{g}$ be the corresponding real Lie algebra.

Recall that $\mathcal{H}$ has an orthogonal basis indexed by a set of paths via the bijection given in Lemma~\ref{pathbij}.

Let $K \subset \U(\mathcal{H})$ be the maximal torus subgroup of diagonal matrices.  Note that the corresponding Cartan subalgebra is a maximal compact Cartan subalgebra.
Using \eqref{eq:twistuniform}, one computes
$$ J_i \cdot {\bf p}=q^{\frac{x_i^2-x_{i-1}^2-\alpha^2+1}{2}}.$$

For $k=n-3$, the vector space is $(n-1)$-dimensional.  Recall the ordered basis $\{e_1, \ldots, e_{n-1} \}$ from earlier where 
\[
e_i=(\alpha, 2\alpha+1, 3\alpha+2, \ldots, i\alpha+(i-1), 
(i+1)\alpha+(i-2), (i+2)\alpha+(i-1), \ldots, n\alpha+(n-3)) \ .
\]
In this basis, recall the matrix for $J_j$ (for $j=1,\ldots, n-1$) from \eqref{eq:JMn-3}.
These $J_i$ generate a dense subgroup of $K$.  
Let $\mathfrak{h}=\Lie(K) \subset \mathfrak{g}$.
Note that $\mathfrak{g} \subset \mathfrak{u}(\mathcal{H})$ has multiplicity-free root spaces.

Any $\mathfrak{g}$-submodule of $\mathcal{H}$ is a weight module.
If $\mathcal{H}$ is an irreducible $\mathfrak{g}$-module, then $\mathfrak{g}=\mathfrak{u}(\mathcal{H})$.

Consider the graph of paths where the vertices are paths and two vertices are connected if the corresponding paths differ by one labeled edge in their 
corresponding trees.
If paths ${\bf p}$ and ${\bf p}'$ are connected by an edge in the graph, then ${\bf p}' \in \mathfrak{g} \cdot {\bf p}$ using operators $\sigma_i$ and $J_i$ for
some $i$.  Since the graph of paths is connected, we obtain $\mathfrak{g}=\mathfrak{u}(\mathcal{H})$.

When $\U(\mathcal{H})$ is a compact Lie group, it is well known that the exponential map is surjective.  This is no longer true in the non-compact case.  For example, it fails for $\SL(2,\R) \cong \SU(1,1)$.  There are suitable variations of this property that suit our purposes.  
The exponential map is surjective for $\PSU(p,p)$ \cite[Main Theorem]{DokNguy}.  See also \cite[Theorem 4.7]{DokHof1} for a survey.  In that theorem, there is a surjectivity result for a certain quotient of $\SU(p,q)$ with $p \neq q$.

The image of the exponential map is actually dense for the adjoint group of $\mathfrak{su}(p,q)$ with $p \neq q$ \cite[Theorem C]{DokNguy2}.   By \cite[Theorem 4.5]{DokNgConj}, the image of the exponential map is actually dense in $\SU(p,q)$ for $p \neq q$.  See also \cite[Section IV]{Neeb} for similar results.



Thus using the commutative diagram \eqref{diag:exp}, $G=\U(\mathcal{H})$ and in particular, $\Br_n$ generates a dense subgroup of $\PSU(\mathcal{H})$.

\begin{equation}\label{diag:exp}
\xy
 (-20,8)*+{G}="tl";
 (20,8)*+{\U(\mathcal{H})}="tr";
 (-20,-8)*+{\mathfrak{g}}="bl";
 (20,-8)*+{\mathfrak{u}(\mathcal{H})}="br";
    {\ar^{\exp} "bl";"tl"};
    {\ar_{\cong} "bl";"br"};
    {\ar_{\exp} "br";"tr"};
    {\ar^{} "tl";"tr"};
\endxy
\end{equation}
\end{proof}

\begin{remark}
Moskowitz proves \cite[Theorem 2.1]{Mosk1} that the exponential map is surjective for any rank 1 centerless, connected, non-compact Lie group.  This was extended to the higher rank case at the expense of additional hypotheses imposed on Cartan subgroups \cite[Theorem 3.2]{Mosk1}.

For results about density of the exponential map for non-semisimple Lie groups, see \cite{HofMuk}.
\end{remark}

\section{Perspectives on quantum computation} \label{sec:quantum}
In this section, we consider the reduced Burau representation for the braid group on three strands and restrict ourselves to values of $\alpha$ that make the Hermitian pairing positive definite.  Thus, the image of the braid group in this section lies in $\PSU(2)$.  This then serves as a model for a single qubit quantum computer.

Consider the 2-dimensional vector space $\mathcal{H}_{3,\alpha}$ spanned by $ \{v_1, v_2 \}$ where
\begin{equation} \label{eq:basis} v_1 =
 \hackcenter{
\begin{tikzpicture}[yscale=-1, scale=0.7,  decoration={markings, mark=at position 0.6 with {\arrow{>}};},]
\draw[very thick, postaction={decorate}] (0,0) to [out=90, in=220] (.75,1);
\draw[very thick, postaction={decorate}] (1.5,-1) to [out=90, in=-30] (.75,1);
\draw[very thick,  postaction={decorate}] (.8,1) to  (.8,2);
\draw[very thick,  postaction={decorate}] (-.5,-1)to [out=90, in=210] (0,0);
\draw[very thick, postaction={decorate}] (.5,-1) to [out=90, in=-30] (0,0);
\node at (.75,-.8) {$\scs \alpha$};
\node at (-.75,-.8) {$\scs \alpha$};
\node at (1.75,-.8) {$\scs \alpha$};
\node at (.5,1.6) {$\scs 3 \alpha$};
\node at (-.5,.6) {$\scs 2 \alpha+1$};
\end{tikzpicture} }
\quad \quad
v_2 =
 \hackcenter{
\begin{tikzpicture}[yscale=-1, scale=0.7,  decoration={markings, mark=at position 0.6 with {\arrow{>}};},]
\draw[very thick, postaction={decorate}] (0,0) to [out=90, in=220] (.75,1);
\draw[very thick, postaction={decorate}] (1.5,-1) to [out=90, in=-30] (.75,1);
\draw[very thick,  postaction={decorate}] (.8,1) to  (.8,2);
\draw[very thick,  postaction={decorate}] (-.5,-1)to [out=90, in=210] (0,0);
\draw[very thick, postaction={decorate}] (.5,-1) to [out=90, in=-30] (0,0);
\node at (.75,-.8) {$\scs \alpha$};
\node at (-.75,-.8) {$\scs \alpha$};
\node at (1.75,-.8) {$\scs \alpha$};
\node at (.5,1.6) {$\scs 3 \alpha$};
\node at (-.5,.6) {$\scs 2 \alpha-1$};
\end{tikzpicture} } \ .
\end{equation}

This is an orthogonal basis with pairing from Proposition \ref{prop:pairinggen}:
\begin{equation} \label{formvi}
n_1^2:=\langle v_1 | v_1 \rangle =
\frac{[2\alpha]}{[\alpha+1]}=2\sin(\frac{\pi\alpha}2) ,
\quad  \quad
n_2^2:=\langle v_2 | v_2 \rangle =
\frac{-[2\alpha] [3\alpha+1]}{[\alpha+1]^2}=2\sin(\frac{\pi\alpha}2)(1-2\cos(\pi\alpha)) \ .
\end{equation}

Define a normalizing matrix to go from the ${v_1}, {v_2}$ basis to an orthonormal basis by:
\[
N=
\begin{pmatrix}
\frac{1}{n_1} & 0 \\
0 & \frac{1}{n_2}
\end{pmatrix} \ .
\]

Consider the following basis of $\mathcal{H}_{3,\alpha}$:
\[ v_1' =
 \hackcenter{
\begin{tikzpicture}[xscale=-1, yscale=-1, scale=0.7,  decoration={markings, mark=at position 0.6 with {\arrow{>}};},]
\draw[very thick, postaction={decorate}] (0,0) to [out=90, in=220] (.75,1);
\draw[very thick, postaction={decorate}] (1.5,-1) to [out=90, in=-30] (.75,1);
\draw[very thick,  postaction={decorate}] (.8,1) to  (.8,2);
\draw[very thick,  postaction={decorate}] (-.5,-1)to [out=90, in=210] (0,0);
\draw[very thick, postaction={decorate}] (.5,-1) to [out=90, in=-30] (0,0);
\node at (.75,-.8) {$\scs \alpha$};
\node at (-.75,-.8) {$\scs \alpha$};
\node at (1.75,-.8) {$\scs \alpha$};
\node at (.5,1.6) {$\scs 3 \alpha$};
\node at (-.5,.6) {$\scs 2 \alpha+1$};
\end{tikzpicture} }
\quad \quad
v_2' =
 \hackcenter{
\begin{tikzpicture}[xscale=-1, yscale=-1, scale=0.7,  decoration={markings, mark=at position 0.6 with {\arrow{>}};},]
\draw[very thick, postaction={decorate}] (0,0) to [out=90, in=220] (.75,1);
\draw[very thick, postaction={decorate}] (1.5,-1) to [out=90, in=-30] (.75,1);
\draw[very thick,  postaction={decorate}] (.8,1) to  (.8,2);
\draw[very thick,  postaction={decorate}] (-.5,-1)to [out=90, in=210] (0,0);
\draw[very thick, postaction={decorate}] (.5,-1) to [out=90, in=-30] (0,0);
\node at (.75,-.8) {$\scs \alpha$};
\node at (-.75,-.8) {$\scs \alpha$};
\node at (1.75,-.8) {$\scs \alpha$};
\node at (.5,1.6) {$\scs 3 \alpha$};
\node at (-.5,.6) {$\scs 2 \alpha-1$};
\end{tikzpicture} }
\]


By Lemma \ref{6jformulalemma},
\[
v_1=- \qd(2\alpha+1) [\alpha+1] v_1' - \qd(2\alpha-1) [\alpha+1] v_2'
\]
\[
v_2=\qd(2\alpha+1) [3\alpha+3] v_1' - \qd(2\alpha-1) [\alpha+1] v_2'
\]
and so the change of basis matrix between the basis $\{v_1, v_2 \}$ and $ \{v_1', v_2' \}$ is given by
\[
A=
\begin{pmatrix}
- \qd(2\alpha+1) [\alpha+1]
& \qd(2\alpha+1) [3\alpha+3]  \\
   - \qd(2\alpha-1) [\alpha+1]
& - \qd(2\alpha-1) [\alpha+1]
\end{pmatrix}
\;\; = \;\;
\frac{1}{[2\alpha]}
\begin{pmatrix}
-    [\alpha+1]
&   [3\alpha+3]   \\
   [\alpha+1]
&     [\alpha+1]
\end{pmatrix}
\ .
\]
Note that $\md(2\alpha+1) = -1/[2\alpha+2] = 1/[2\alpha]$ and $\md(2\alpha-1)=-1/[2\alpha]$.


\begin{lemma} \label{lem:goodalpha}
Assume $\alpha$ is irrational.
The vector space $ \mathcal{H}_{3,\alpha}$ has a definite Hermitian inner product if and only if $ \frac{1}{3} +2k < \alpha < \frac{5}{3}+2k$ for an integer $k$.
Furthermore, this form is positive definite if and only if 
$ \frac{1}{3} +4k < \alpha < \frac{5}{3}+4k$ for an integer $k$.
\end{lemma}

\begin{proof}
This is a straightforward computation using \eqref{formvi}.
\end{proof}

From Section \ref{sec:gen}, there is an action of the braid group $\Br_3$ where $\sigma_1$ crosses the two top left strands and $\sigma_2$ crosses the two top right strands.
Just as in the previous section, after a lengthy calculation, using Lemma \ref{Br2lemma} and a change of basis from $\{{v_1},{v_2} \} $ to $\{{v_1'},{v_2'} \} $ one gets the matrices for $\sigma_1$ and $\sigma_2 $ in the basis $\{{v_1}, {v_2}\}$ to be
\[
{\sigma_1}=
\begin{pmatrix}
q^{\frac{(\alpha+1)^2}{2}} & 0 \\
0 & q^{\frac{(\alpha-1)^2}{2}}
\end{pmatrix} ,
\]


\begin{align*}
{\sigma_2}&=A^{-1} \sigma_1 A \\
&=
\frac{1}{[\alpha+1]+[3\alpha+3]}
\begin{pmatrix}
[3\alpha+3] q^{\frac{(\alpha-1)^2}{2}}  + [\alpha+1]  q^{\frac{(\alpha+1)^2}{2}}
& & [3 \alpha+3] q^{\frac{(\alpha-1)^2}{2}}  - [3 \alpha+3]  q^{\frac{(\alpha+1)^2}{2}}  \\
[\alpha+1] q^{\frac{(\alpha-1)^2}{2}}  - [\alpha+1]  q^{\frac{(\alpha+1)^2}{2}}
& & [\alpha+1] q^{\frac{(\alpha-1)^2}{2}} + [3 \alpha+3] q^{\frac{(\alpha+1)^2}{2}}
\end{pmatrix} .
\end{align*}



\begin{remark}
The matrices $ N^{-1} \sigma_1 N^{}$ and $ N^{-1} \sigma_2 N^{}$ are matrices for the braiding operators in the orthonormal basis.
\end{remark}

In the Fibonacci category, the half-twist braiding operator $\sigma_1 \sigma_2 \sigma_1$ gives rise to the Hadamard gate in $\PSU(2)$. 
We record here a formula for this braid in our category of interest.  Note that when $\alpha=\frac{1}{2}$, the half-twist is the Hadamard operator up to a scalar.  However, we exclude such an $\alpha$ from this section since it is rational. 

\begin{lemma}
The braid group generators of $\Br_3$ associated to the parameter $\alpha$ satisfy
\begin{align} \label{eq:halftwist}
\sigma_1\sigma_2\sigma_1 &=
\frac{(-1)^{3/4} e^{\frac{1}{4} i \pi  \alpha  (3 \alpha +2)}}{1-e^{i \pi  \alpha }}
\left(
\begin{array}{cc}
 1 & -1+2\cos(\pi \alpha)  \\
 -1 & -1 \\
\end{array}
 \right) \ ,
 \\  \label{s1s2s22}
(\sigma_1\sigma_2\sigma_1)^2 &=  i e^{\frac{3}{2} i \pi  \alpha ^2}\left(
\begin{array}{cc}
 1& 0 \\
 0 & 1 \\
\end{array}
\right) .
\end{align} 
Thus $\sigma_1 \sigma_2 \sigma_1$ is an involution in $\PSU(2)$.

In the orthonormal basis, these operators have the form:
\begin{equation} \label{eq:normdefHad}
N^{-1} \sigma_1\sigma_2\sigma_1 N=
 \frac{(-1)^{3/4} e^{\frac{1}{4} i \pi  \alpha  (3 \alpha +2)}}{1-e^{i \pi  \alpha }}
\left(
\begin{array}{cc}
 1 & \varsigma \\
 \varsigma & -1 \\
\end{array}
\right), \qquad \varsigma= -\frac{n_2}{n_1}\ , \\
\end{equation}
  with $n_1$ and $n_2$ defined in \eqref{formvi}, so that  $\varsigma^2=1-2\cos(\pi\alpha)$ and \eqref{s1s2s22} becomes 
\begin{equation}
N^{-1} (\sigma_1\sigma_2\sigma_1)^2 N =\left(
\begin{array}{cc}
 i e^{\frac{3}{2} i \pi  \alpha ^2} & 0 \\
 0 & i e^{\frac{3}{2} i \pi  \alpha ^2} \\
\end{array}
\right) \ .
\end{equation}

\end{lemma}

Next, we consider the following bi-invariant metric on $\PSU(2)$:
\begin{equation} \label{eq:metric}
    d(x,y)= \sqrt{1-\frac{|\Tr(x^\dagger y)|}{2}} \ .
\end{equation}
Earlier it was shown that $\sigma_1, \sigma_2$ generate a dense subgroup of $\PSU(2)$.  Then it is easy to see that the operators $\sigma_1 $ and $N^{-1} \sigma_1 \sigma_2 \sigma_1 N$ from \eqref{eq:normdefHad} also generate a dense subgroup.

The Solovay-Kitaev theorem provides an algorithm to approximate up to an error $\epsilon$ a given element $M \in \PSU(2)$ by elements in our generating set $\sigma_1, N^{-1} \sigma_1 \sigma_2 \sigma_1 N$ (and their inverses), using a word of order $\Log^c(\frac{1}{\epsilon})$ where $c \geq 1$.  It is well known that $c$ cannot be less than $1$.  We believe that with our generating set, we could obtain the optimal bound of $c=1$.



\begin{figure}[htbp]
    \centering
    \begin{subfigure}[b]{0.48\textwidth}
        \centering
        \includegraphics[width=\textwidth]{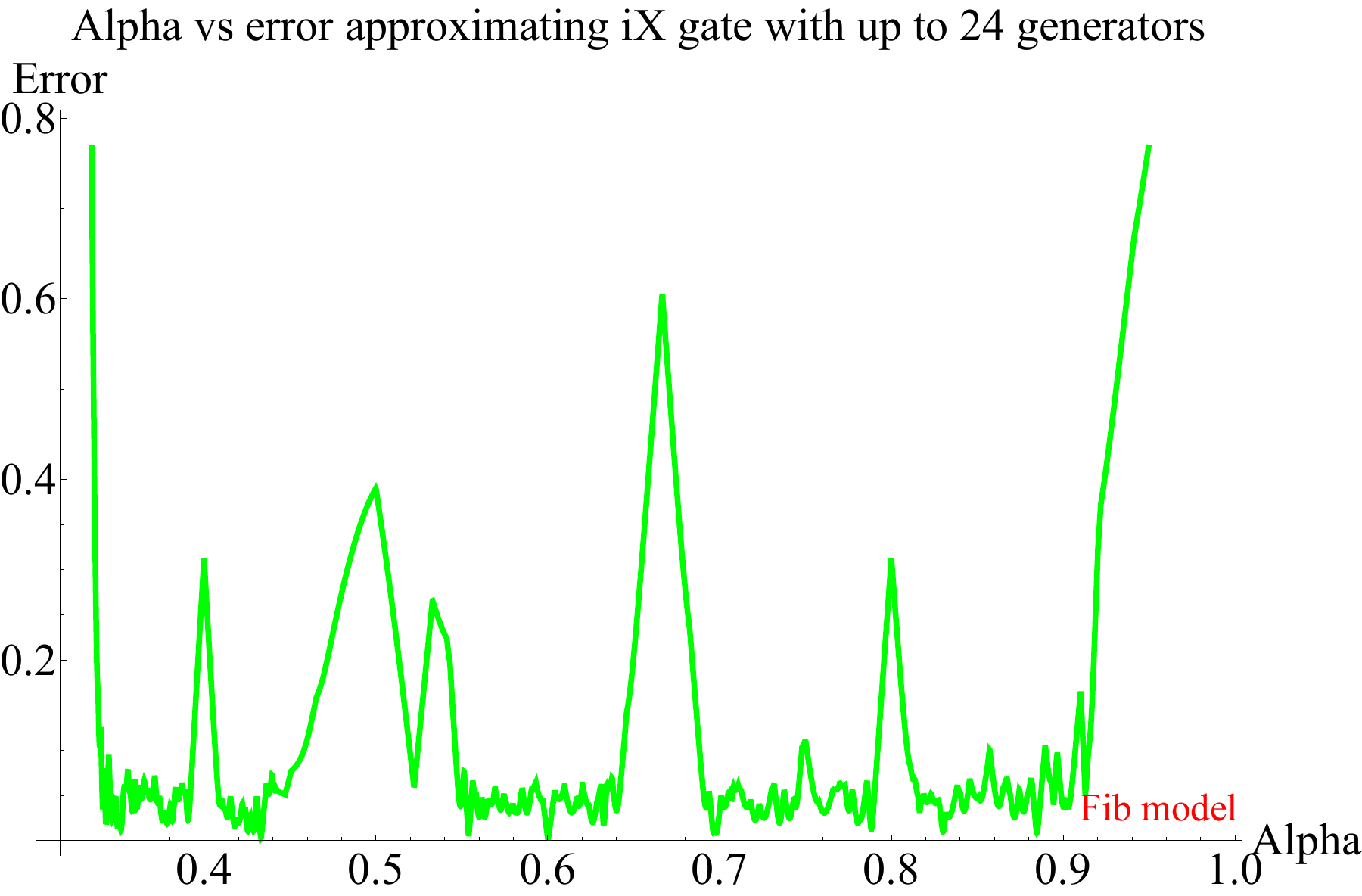}
        \caption{Minimum for Fibonacci model is $3\times 10^{-3}$ versus $7 \times 10^{-4}$ occurring at $\alpha=0.6002$. Most values of $\alpha$ underperform the Fibonacci model.}
        \label{fig:iX}
    \end{subfigure}
    \hfill 
    \begin{subfigure}[b]{0.48\textwidth}
        \centering
        \includegraphics[width=\textwidth]{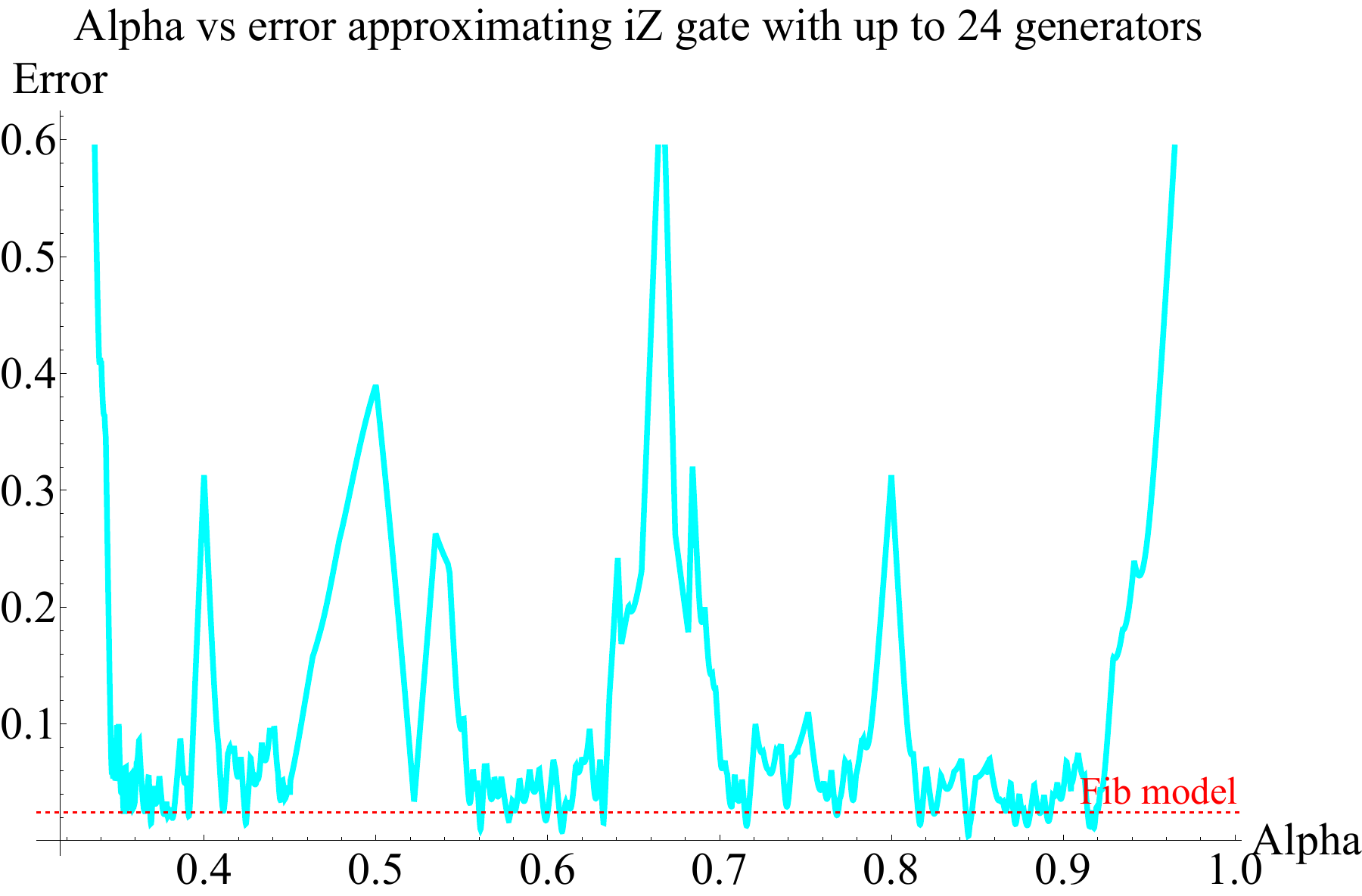}
        \caption{Minimum for Fibonacci model is $2 \times 10^{-2}$ versus $4 \times 10^{-3}$ occurring at $\alpha =0.8448$.}
        \label{fig:iZ}
    \end{subfigure}
    \vspace{1em} 

    \begin{subfigure}[b]{0.48\textwidth}
        \centering
        \includegraphics[width=\textwidth]{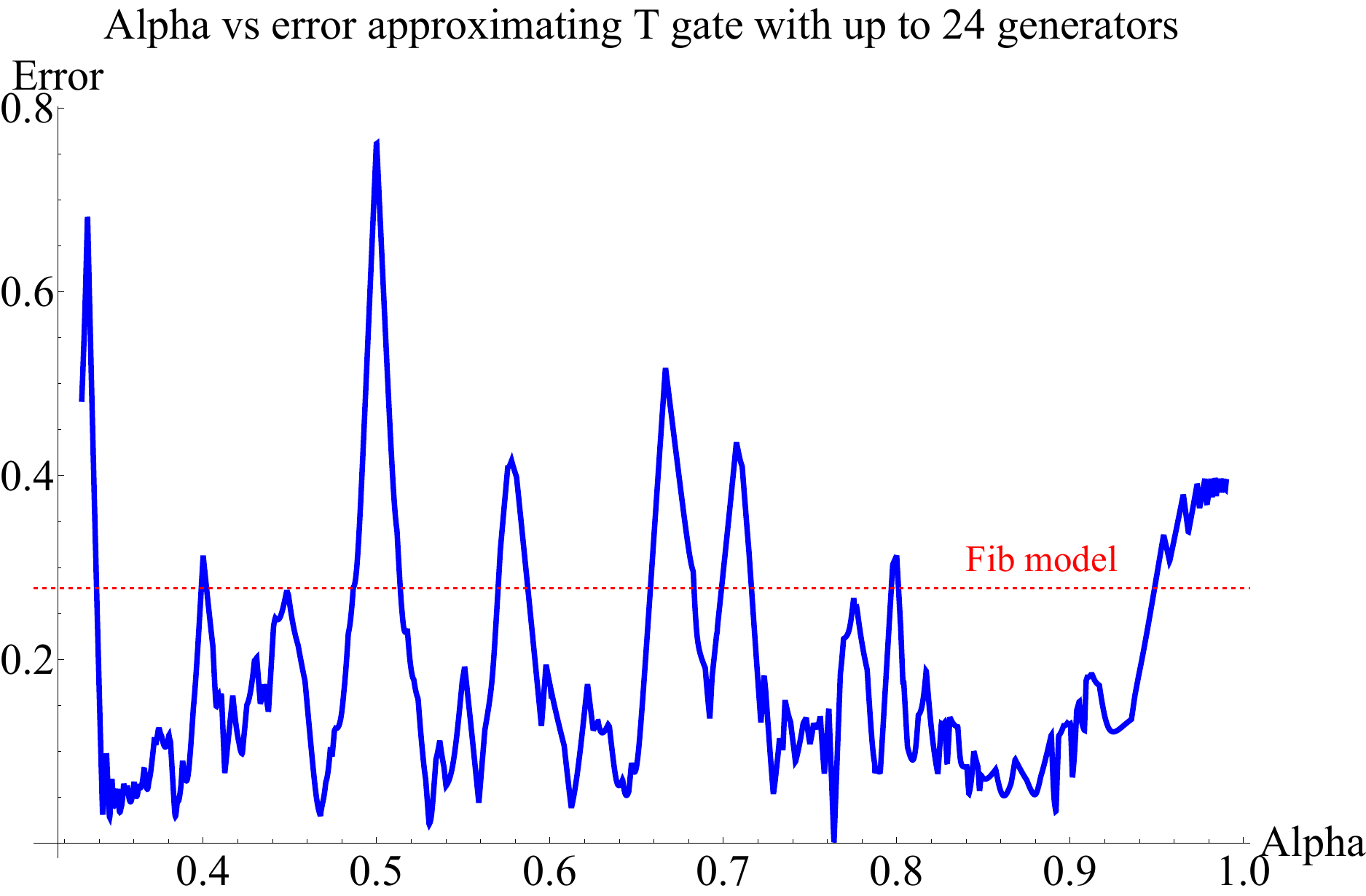}
        \caption{Minimum for Fibonacci model is $2 \times 10^{-1}$ versus $3 \times 10^{-3}$ for $\alpha=0.7639$.}
        \label{fig:T}
    \end{subfigure}
    \hfill
    \begin{subfigure}[b]{0.48\textwidth}
        \centering
        \includegraphics[width=\textwidth]{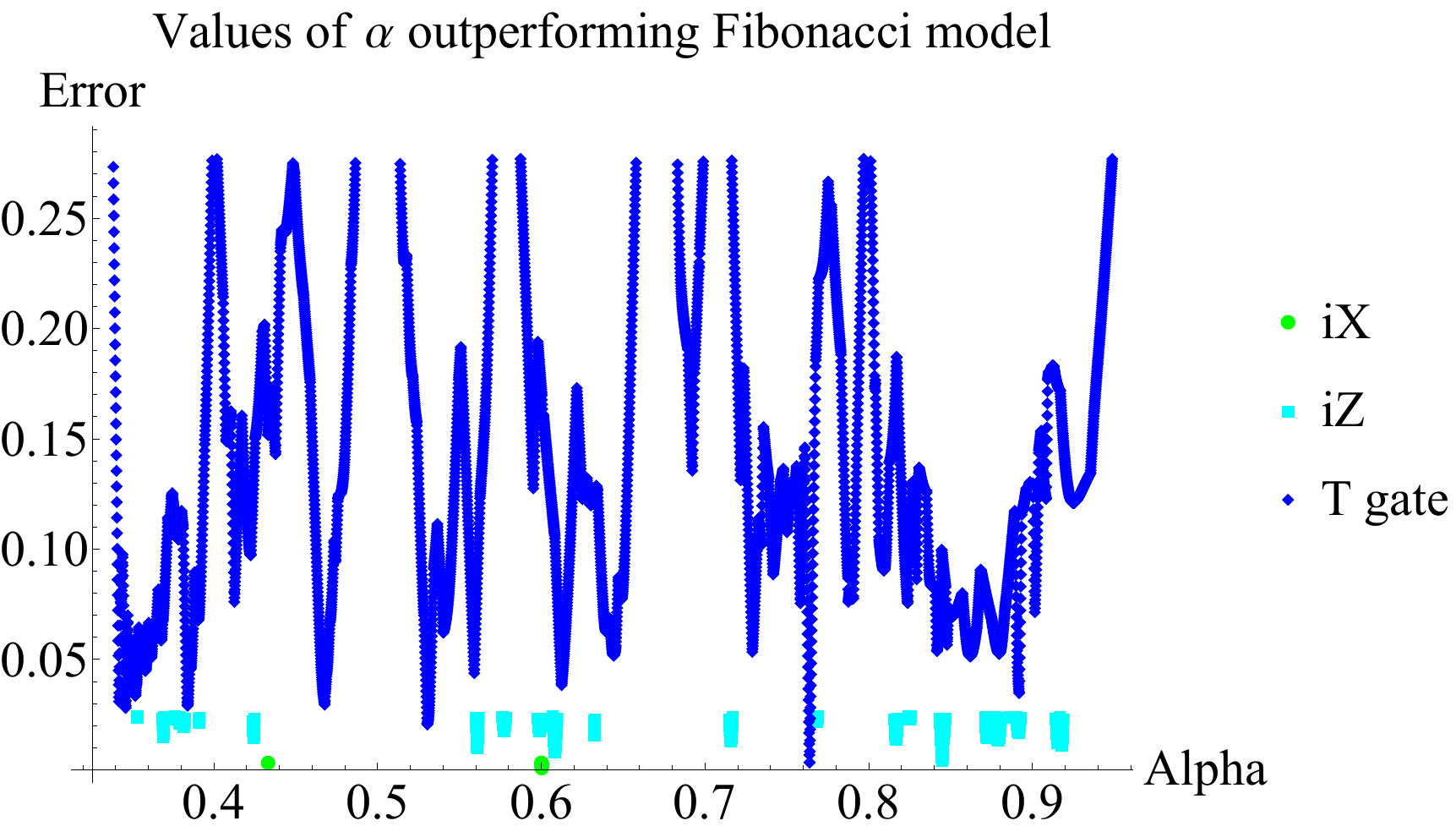}
        \caption{Values of $\alpha$ where the weaves of up to 24 generators on the basis \eqref{eq:basis} outperform those of the Fibonacci model. Values near $\alpha=0.6$ outperform the Fibonacci model by an order of magnitude for all gates tested.  }
        \label{fig:Mix}
    \end{subfigure}

    \caption{
    %
 For each fixed value of $\alpha$, we perform a brute-force search over all weaves of length up to 24 to find the unitary transformation that best approximates a given target gate, measured using the operator norm. The plot shows the minimal error $\epsilon(U(\mathbf{n}), T)$ achieved for unitaries $U(\mathbf{n})$ constructed from weaves as in equation (5.8) with target gate $T$, across the sampled values of $\alpha$.
    }
    \label{fig:all-plots}
\end{figure}

In Figure~\ref{fig:all-plots}, we give numerical evidence illustrating the efficiency of approximating elements of $\PSU(2)$ using braiding operators on the basis from \eqref{eq:basis}. 

We consider a brute force search over the best possible approximation achievable using a fixed number $N$ of braid generators.  
The quality of the approximation is measured using the operator norm where the distance between operators $U$ and $V$ is given by $\epsilon(U,V)=\|U-V\|$ where $\|O\|$ is the square root of highest eigenvalue of $O^{\dagger}O$.   
In this example, we take $N=24$ and consider approximations of the $\PSU(2)$ matrices $iX$, $iZ$, and $T$ (or $\pi/8$ phase).   We perform such a brute force search for varying values of $\alpha$ and plot the minimal error achieved with the fixed $N=24$ braid generators.  

We follow \cite{PhysRevB.75.165310} and restrict our search to a subclass of all possible braids known as weaves.  Weaves are braids in which a single mobile quasiparticle interlaces between fixed and immutable quasiparticles. It was argued in  \cite{PhysRevLett.96.070503} that this restricted class of braids could still be used for universal quantum computation and may have advantages in physical implementations.   
Hence, we consider approximations of a target unitary by unitaries $ U(\mathbf{n} )$ obtained from weaves taking the form 
\begin{equation} \label{eq:weave}
  U(\mathbf{n} ):=\sigma_1^{n_m}\sigma_2^{n_{m-1}} \dots \sigma_1^{n_3}\sigma_2^{n_2} \sigma_1^{n_1},
\end{equation}
for $\mathbf{n}=(n_m,n_{m-1},\dots, n_2,n_1)$, where
 we can assume that all the $n_i$ are even integers and that $n_i \neq 0$ for $1< i <m$.  The exponents being even translates into the mobile quasiparticle wrapping entirely around one of the neighboring quasiparticles and returning to the center position.  

As a point of reference, we have also included the best approximation obtainable from braiding on a similarly defined fusion basis in the context of Fibonacci anyons associated with $\SU(2)$ at a 5th root of unity.  This is a popular model for universal topological quantum computation.     For all of the gates considered, the braiding defined on the fusion basis \eqref{eq:basis} can achieve an order of magnitude improvement over the best approximation obtained in the Fibonacci model for some value of $\alpha$.  

The Fibonacci result is depicted as a horizontal line in the graphs representing the best achievable error by weaves with 24 generators in the Fibonacci model.  In Figure~\ref{fig:iX}, there are only a few values of $\alpha$ which improve upon the Fibonacci model's best approximation.   However, in Figure~\ref{fig:iZ}, there is a much larger range of values of $\alpha$ that produce improvements over the Fibonacci model in approximating the $iZ$ gate.  Figure~\ref{fig:T} shows that most values of $\alpha$ give a much better approximation of the  $\pi/8$ phase gate than the Fibonacci model.  Finally, in Figure~\ref{fig:Mix}, we graph all the values of $\alpha$, which offer improvements over the Fibonacci model.   The values around $\alpha =0.6$ offer improvements over the Fibonacci model for all gates considered.

\begin{remark}
Finding an optimal set of topological generators is a challenging question that has connections to many different areas of mathematics.  There are even many different notions of optimal.
One particular perspective is the notion of efficiently universal mentioned earlier, meaning that it saturates the lower bound of the exponent in the Kitaev-Solovay Theorem.  In the work of Chuang, Harrow, and Recht \cite{HRC}, it was shown that a finite topologically generating set saturates this bound if their corresponding mixing operator on $L^2(\SU(2))$ is close in some sense to the projection operator onto constant functions.  Using the work of Lubotzky, Phillips, and Sarnak \cite{LPS1, LPS2}, they deduce that the so-called $V$-gates \eqref{eq:Vgates} 
\begin{equation} \label{eq:Vgates}
V_x= \frac{1}{\sqrt{5}}
\begin{pmatrix}
1+2i & 0 \\
0 & 1-2i
\end{pmatrix},
\quad \quad
V_y= \frac{1}{\sqrt{5}}
\begin{pmatrix}
1 & 2 \\
-2 & 1
\end{pmatrix},
\quad \quad
V_z= \frac{1}{\sqrt{5}}
\begin{pmatrix}
1 & 2i \\
2i & 1
\end{pmatrix}
\end{equation}
are efficiently universal. For explicit algorithms, see \cite{BocGurSv, Ross, BlassBochGur}.
These algorithms use solutions of certain grid problems and diophantine equations.

Before explicit algorithms were constructed for these $V$-gates, analogous questions for the more standard Clifford $+T$-gates 
\begin{equation} \label{eq:cliffgates}
H= \frac{1}{\sqrt{2}}
\begin{pmatrix}
1 & 1 \\
1 & -1
\end{pmatrix},
\quad \quad
T= 
\begin{pmatrix}
1 & 0 \\
0 & e^{\pi i / 4}
\end{pmatrix}  ,
\end{equation}
were studied in \cite{KMM1, GS1, KMM2, Sel1, KMM3, RS}.

For the Fibonacci category, it was shown in \cite{KBS} that the generating gate set coming from the standard $\sigma_1, \sigma_2$ (or an equivalent set) is efficiently universal, in the sense that in order to approximate an arbitrary element of $\PSU(2)$ up to an error $\epsilon$, one needs a braid whose length is order $\Log(\frac{1}{\epsilon})$.  This is the fastest possible speed-up that could be achieved from the Solovay-Kitaev theorem.  Their proof is number-theoretic.  

Based on experimental evidence, we conjecture that our operators $$\sigma_1, \sigma_1^{-1}, N^{-1} \sigma_1 \sigma_2 \sigma_1 N, N^{-1} \sigma_1^{-1} \sigma_2^{-1} \sigma_1^{-1} N$$ are efficiently universal.
It is well known that an optimal generator of $\U(1)$ is $e^{2 \pi i \phi}$ where $\phi$ is the golden ratio (up to an integer).
See \cite{GrahamLint, Stier1} for a more in-depth analysis.
This then tells us which $\alpha=4\phi$ modulo $1$ gives the best approximations of diagonal matrices, which in turn may improve the length of our quantum circuit even more.

There are also notions of golden and super golden gate sets due to Sarnak \cite{Sarletter} and developed in \cite{EP,PS-golden}.  By definition, these gates cover the Lie group in an optimal way and there are efficient algorithms to write an arbitrary element in the Lie group in terms of the generating gate set.  For a precise definition see \cite[Definition 2.8]{EP}.  
Proving that a gate set is golden or super golden uses some deep theorems in number theory.  The so-called icosahedral super golden gates enjoy certain advantages over other sets \cite{BlackSt}.
It would be interesting to prove that the gates we consider here are golden.  They cannot be super golden since, by definition, the elements of the generating gate set must be of finite order, and the elements $\sigma_1, \sigma_2$ that we consider certainly do not possess this property.
\end{remark}

\section{Unitary representations of the braid group: singular part of the category} \label{sec:sing}
In this section we will consider the singular part of the category from Section~\ref{SS:UqH}.  Representations of the braid group coming from morphism spaces here, lie in discrete subgroups of the corresponding unitary groups and thus have no value for quantum computation.

The simple module $V_0$ in the singular part of the category is projective.
Recall that in Figure \ref{fig:P_i}, a basis for the module $P_0$ which is the projective cover of $\C^H_0$ is given.
There is a unique isomorphism up to scalar $\phi \colon V_0 \otimes V_0 \cong P_0$ where
\[
\phi(v_0 \otimes v_0) = w^R, \quad
\phi(v_1 \otimes v_1) = w^L, \quad
\phi(v_0 \otimes v_1) = w^H, \quad
\phi(v_1 \otimes v_0) = w^S+q w^H .
\]
There is an isomorphism $V_0^{\otimes 3} \cong V_0 \otimes P_0 \cong (\C_2 \oplus \C_0^{\oplus 2} \oplus \C_{-2}) \otimes V_0$.

So, we have a potential model for a qubit
\[
\Hom(V_0^{}, V_0^{\otimes 3}) = \Span \{ {\pi_1}, {\pi_2} \}
\]
where $\pi_1$ and $\pi_2$ are defined as follows:
\[
\pi_1^{}(v_0)=  v_1 \otimes v_0 \otimes v_0 + q^{-1} v_0 \otimes v_1 \otimes v_0, \quad \quad
\pi_1^{}(v_1)= v_1 \otimes v_0 \otimes v_1 + q^{-1} v_0 \otimes v_1 \otimes v_1  ,
\]
\[
\pi_2^{}(v_0)= v_0 \otimes v_1 \otimes v_0 + q^{-1} v_0 \otimes v_0 \otimes v_1, \quad \quad
\pi_2^{}(v_1)= v_1 \otimes v_1 \otimes v_0 + q^{-1} v_1 \otimes v_0 \otimes v_1.
\]
We depict these morphisms by:
\[ {\pi_1} =
 \hackcenter{
\begin{tikzpicture}[yscale=-1, scale=0.7,  decoration={markings, mark=at position 0.6 with {\arrow{>}};},]
\draw[very thick, postaction={decorate}] (0,0) to [out=90, in=220] (.75,1);
\draw[very thick, postaction={decorate}] (1.5,-1) to [out=90, in=-30] (.75,1);
\draw[very thick,  postaction={decorate}] (.8,1) to  (.8,2.3);
\draw[very thick,  postaction={decorate}] (-.5,-1)to [out=90, in=210] (0,0);
\draw[very thick, postaction={decorate}] (.5,-1) to [out=90, in=-30] (0,0);
\node at (.75,-.8) {$\scs V_0$};
\node at (-.75,-.8) {$\scs V_0$};
\node at (1.75,-.8) {$\scs V_0$};
\node at (.4,2) {$\scs V_0$};
\node at (-.5,.6) {$\scs P_0$};
\node[draw, thick, fill=blue!20,rounded corners=4pt,inner sep=4pt] (X) at (.75,1) {$1$};
\end{tikzpicture} }
\quad \quad
{\pi_2} =
 \hackcenter{
\begin{tikzpicture}[yscale=-1, scale=0.7,  decoration={markings, mark=at position 0.6 with {\arrow{>}};},]
\draw[very thick, postaction={decorate}] (0,0) to [out=90, in=220] (.75,1);
\draw[very thick, postaction={decorate}] (1.5,-1) to [out=90, in=-30] (.75,1);
\draw[very thick,  postaction={decorate}] (.8,1) to  (.8,2.3);
\draw[very thick,  postaction={decorate}] (-.5,-1)to [out=90, in=210] (0,0);
\draw[very thick, postaction={decorate}] (.5,-1) to [out=90, in=-30] (0,0);
\node at (.75,-.8) {$\scs V_0$};
\node at (-.75,-.8) {$\scs V_0$};
\node at (1.75,-.8) {$\scs V_0$};
\node at (.4,2) {$\scs V_0$};
\node at (-.5,.6) {$\scs P_0$};
\node[draw, thick, fill=blue!20,rounded corners=4pt,inner sep=4pt] (X) at (.75,1) {$2$};
\end{tikzpicture} } \ ,
\]
where the vertices labeled $1$ and $2$ denote the induced maps $P_0 \otimes V_0 \rightarrow V_0$ coming from the maps $\pi_1$ and $\pi_2$ under the identification $\phi \colon V_0 \otimes V_0 \cong P_0$.

In this basis, the matrices of the braid group generators $\sigma_1$ and $\sigma_2$ are given by:
\begin{equation} \label{eq:sigmassing}
{\sigma_1}=
q^{\frac{1}{2}}
\begin{pmatrix}
1 & -q \\
0 & 1
\end{pmatrix} ,
\quad \quad
\sigma_2=
q^{\frac{1}{2}}
\begin{pmatrix}
1 & 0 \\
-q & 1
\end{pmatrix} .
\end{equation}
\begin{lemma}
The matrices $\sigma_1$ and $\sigma_2$ have infinite order.
\end{lemma}

\begin{proof}
This is clear from inspection.
\end{proof}

Then one could compute their Hermitian adjoints in $\Hom(V_0^{\otimes 3}, V_0^{})$
\begin{align*}
{\pi_1^{\dagger}}(v_0 \otimes v_0 \otimes v_0) = 0 & & {\pi_2^{\dagger}}(v_0 \otimes v_0 \otimes v_0) = 0 \\
{\pi_1^{\dagger}}(v_1 \otimes v_0 \otimes v_0) = q v_0 & &{\pi_2^{\dagger}}(v_1 \otimes v_0 \otimes v_0) = 0 \\
{\pi_1^{\dagger}}(v_0 \otimes v_1 \otimes v_0) = v_0 & &{\pi_2^{\dagger}}(v_0 \otimes v_1 \otimes v_0) = qv_0 \\
{\pi_1^{\dagger}}(v_0 \otimes v_0 \otimes v_1) = 0 & & {\pi_2^{\dagger}}(v_0 \otimes v_0 \otimes v_1) = v_0 \\
{\pi_1^{\dagger}}(v_1 \otimes v_1 \otimes v_0) = 0 & & {\pi_2^{\dagger}}(v_1 \otimes v_1 \otimes v_0) = q v_1 \\
{\pi_1^{\dagger}}(v_1 \otimes v_0 \otimes v_1) = q v_1 & & {\pi_2^{\dagger}}(v_1 \otimes v_0 \otimes v_1) = v_1 \\
{\pi_1^{\dagger}}(v_0 \otimes v_1 \otimes v_1) = v_1 & & {\pi_2^{\dagger}}(v_0 \otimes v_1 \otimes v_1) = 0 \\
{\pi_1^{\dagger}}(v_1 \otimes v_1 \otimes v_1) = 0 & & {\pi_2^{\dagger}}(v_1 \otimes v_1 \otimes v_1) = 0 .
\end{align*}

It is straightforward to check
\[\pi_i^{\dagger} \pi_j^{} =
\begin{cases} 
\Id  &\text{ if } i \neq j \\
0 &\text{ otherwise } .
\end{cases}
\]
Thus the vector space $ \Hom(V_0^{}, V_0^{\otimes 3}) $ has a non-degenerate Hermitian pairing
\begin{equation} \label{eq:B}
B = \qd(V_{0})
\begin{pmatrix}
0 & 1 \\
1 & 0
\end{pmatrix} .
\end{equation}

The dual space is:
\[
\Hom(V_0^{\otimes 3}, V_0^{}) = \Span \{ \pi_1^{\dagger}, \pi_2^{\dagger} \}  .
\]
These basis vectors are graphically depicted as follows:
\[ {\pi_1^{\dagger}} =
 \hackcenter{
\begin{tikzpicture}[yscale=-1, scale=0.7]
\draw[very thick, rdirected=.65] (0,-2)  to [out=-90, in=190]  (.75,-3) ;
\draw[very thick, rdirected=.4]  (1.5,-1) to [out=-90, in=-10]  (.75,-3) ;
\draw[very thick, rdirected=.6] (.8,-3) to  (.8,-4);
\draw[very thick,  rdirected=.35] (-.5,-1)to [out=90, in=210] (0,-2);
\draw[very thick, rdirected=.35] (.5,-1) to [out=90, in=-30] (0,-2);
\node at (.75,-.8) {$\scs V_0$};
\node at (-.75,-.8) {$\scs V_0$};
\node at (1.75,-.8) {$\scs V_0$};
\node at (.4,-4) {$\scs V_0$};
\node at (-.5,-2.6) {$\scs P_0$};
\node[draw, thick, fill=blue!20,rounded corners=4pt,inner sep=4pt] (X) at (0,-2) {$1^{\dagger}$};
\end{tikzpicture} } \ ,
\quad \quad
{\pi_2^{\dagger}} =
 \hackcenter{
\begin{tikzpicture}[yscale=-1, scale=0.7]
\draw[very thick, rdirected=.65] (0,-2)  to [out=-90, in=190]  (.75,-3) ;
\draw[very thick, rdirected=.4]  (1.5,-1) to [out=-90, in=-10]  (.75,-3) ;
\draw[very thick, rdirected=.6] (.8,-3) to  (.8,-4);
\draw[very thick,  rdirected=.35] (-.5,-1)to [out=90, in=210] (0,-2);
\draw[very thick, rdirected=.35] (.5,-1) to [out=90, in=-30] (0,-2);
\node at (.75,-.8) {$\scs V_0$};
\node at (-.75,-.8) {$\scs V_0$};
\node at (1.75,-.8) {$\scs V_0$};
\node at (.4,-4) {$\scs V_0$};
\node at (-.5,-2.6) {$\scs P_0$};
\node[draw, thick, fill=blue!20,rounded corners=4pt,inner sep=4pt] (X) at (0,-2) {$2^{\dagger}$};
\end{tikzpicture} } \ .
\]

Recall the endomorphism $x \colon P_0 \rightarrow P_0$ \eqref{def:xmap} which becomes identified with an endomorphism $x \colon V_0 \otimes V_0 \rightarrow V_0 \otimes V_0$.

\begin{lemma}
We have equalities of morphisms
\[
(x \otimes \Id) \circ \pi_2^{} = \pi_1^{} \quad \quad
(x \otimes \Id) \circ \pi_1^{} =0 .
\]
\end{lemma}

It is straightforward to check that the endomorphism $x$ of $V_0 \otimes V_0$ satisfies the property: $x^{\dagger}=x$.

\begin{lemma}
We have equalities of morphisms
\[
\pi_2^{\dagger} \circ (x \otimes \Id) = \pi_1^{\dagger} \quad \quad
\pi_1^{\dagger} \circ (x \otimes \Id) = 0 .
\]
\end{lemma}


One could calculate the matrices for $\sigma_1^{\dagger}, \sigma_2^{\dagger}$ with respect to the indefinite pairing \eqref{eq:B}:
\[
\sigma_1^{\dagger} =
q^{-\frac{1}{2}}
\begin{pmatrix}
1 & q \\
0 & 1
\end{pmatrix} ,
\quad \quad
\sigma_2^{\dagger} =
q^{-\frac{1}{2}}
\begin{pmatrix}
1 & 0 \\
q & 1
\end{pmatrix} .
\]

\begin{lemma}
For $i=1,2$, we have $\sigma_i^{\dagger}=\sigma_i^{-1}$ so the representation of the braid group in $\mathcal{H}_{1,0}$ is unitary.
That is, the image lies in $\PSU(1,1)$.
\end{lemma}

\begin{proof}
This follows immediately from the computations above.
\end{proof}

\begin{remark}
The group $\PSU(1,1)$ is isomorphic to $\PSL(2,\mathbb{R})$.  The image of the braid group is clearly a subgroup isomorphic to $\PSL(2,\mathbb{Z})$ and thus a discrete subgroup of $\PSL(2,\mathbb{R})$.
This lack of density is a further obstruction to doing quantum computing using these representations of the unrolled quantum group.
\end{remark}

More generally, consider the vector space
\[
\mathcal{H}_{n,0} = \Hom(V_0, V_0^{\otimes 2n+1}) .
\]

\begin{theorem} \label{thm:sing}
There is an indefinite unitary representation of the braid group on $\mathcal{H}_{n,0}$.  Furthermore,
$ \dim \mathcal{H}_{n,0} = \binom{2n}{n}$
and has a mixed signature.
\end{theorem}

\begin{proof}
The fact that there is a (possibly indefinite) unitary representation of the braid group follows from \cite[Proposition 5.8]{GLPMS}.
The dimension is determined by a straightforward character calculation.

Part of a basis of the space is given by trees in \eqref{basisHn0}.
\begin{equation} \label{basisHn0}
 \hackcenter{
\begin{tikzpicture}[yscale=-1, scale=0.7,  decoration={markings, mark=at position 0.6 with {\arrow{>}};},]
\draw[very thick, postaction={decorate}] (0,0) to [out=90, in=220] (.75,1);
\draw[very thick, postaction={decorate}] (1.5,-1) to [out=90, in=-30] (.75,1);
\draw[very thick,  postaction={decorate}] (.8,1) to [out=90, in=220]  (1.35,2);
\draw[very thick,  postaction={decorate}] (-.5,-1)to [out=90, in=210] (0,0);
\draw[very thick, postaction={decorate}] (.5,-1) to [out=90, in=-30] (0,0);
\draw[very thick, postaction={decorate}] (2.5,-1) to [out=90, in=-30] (1.35,2);
\node at (1.35,2.2) {$\cdots$};
\draw[very thick, postaction={decorate}] (1.35,2.4) to [out=90, in=220] (1.85,3);
\draw[very thick, postaction={decorate}] (3.5,-1) to [out=90, in=-30] (1.85,3);
\draw[very thick, postaction={decorate}] (1.85,3) to [out=90, in=220] (2.35,4);
\node at (.75,-.8) {$\scs V_0$};
\node at (-.75,-.8) {$\scs V_0$};
\node at (1.75,-.8) {$\scs V_0$};
\node at (.6,2) {$\scs V_0$};
\node at (-.2,.6) {$\scs P_0$};
\node at (3,-.8) {$\cdots$};
\node at (2.25,-.8) {$\scs V_0$};
\node at (3.75,-.8) {$\scs V_0$};
\node at (1,2.6) {$\scs P_0$};
\node at (2,4) {$\scs V_0$};
\node[draw, thick, fill=blue!20,rounded corners=4pt,inner sep=4pt] (X) at (.75,1) {};
\node[draw, thick, fill=blue!20,rounded corners=4pt,inner sep=4pt] (X) at (1.85,3) {};
\end{tikzpicture} }
\end{equation}
At each vertex where two edges labeled $V_0$ meet to produce an edge labeled $P_0$, we assume this vertex is the unique isomorphism (up to scalar) $\phi \colon V_0 \otimes V_0 \rightarrow P_0$.
At each vertex where edges labeled $P_0$ and $V_0$ meet and output $V_0$, there are two linearly independent morphisms.
Thus, this subspace has dimension $2^n$.
The signature of this subspace is calculated easily using the computation for $\mathcal{H}_{1,0}$.
\end{proof}

 



\appendix

\section{Alternative approaches to density} \label{app:alt}
In this section, we give an alternate proof of density for the representation $\mathcal{H}_{n,k,\alpha}$ of the braid group $\Br_n$ for $k=\pm (n-3)$, for $\alpha$ making the Hermitian pairing positive definite.
We begin with an analysis for $n=3$.

Recall the following from~\cite{MR1611329} (see also \cite[Lemma 1]{Cui-Wang}).

\begin{lemma} \label{denseSU2lemma}
Let $U_1$ and $U_2$ are two non-commuting matrices in $\SU(2)$.  If at least one of them has infinite order, then the subgroup generated by $U_1,U_2$ is dense in $\SU(2)$.
\end{lemma}


\begin{lemma} \label{inforderlemma}
The matrices $\sigma_1$ and $\sigma_2$ have infinite order.
\end{lemma}

\begin{proof}
It is clear that for generic real $\alpha$ the matrix $\sigma_1$ has infinite order.  Since $\sigma_2$ is conjugate to $\sigma_1$ ($\sigma_2 = \sigma_1^{-1} \sigma_2^{-1} \sigma_1 \sigma_2 \sigma_1 $) in the braid group, it then follows that $\sigma_2$ also has infinite order.
\end{proof}

\begin{lemma}
The group generated by $\sigma_1$ and $\sigma_2$ is a subgroup of $\PSU(1,1)$ or $\PSU(2)$.
\end{lemma}

\begin{proof}
This is a direct consequence of \cite[Proposition 5.8]{GLPMS}.
\end{proof}

\begin{theorem} \label{mainthm2}
The image of the braid group $\Br_n$ in $\PSU(\mathcal{H}_{n,k,\alpha})$, for $k= \pm (n-3)$ is dense, where $\mathcal{H}_{n, k, \alpha} $ is a non-degenerate finite-dimensional Hermitian vector space which is positive definite.
\end{theorem}

\begin{proof}
We mimic the proof of \cite[Theorems 7.2]{AAEL} in our context of the Burau representation for $k=n-3$.  
Recall the basis $e_1, \ldots, e_{n-1}$ of this space.
The above analysis shows that we have density in the projective unitary group for the space spanned by $e_1, e_2$ using $\sigma_1$ and $\sigma_2$.  In the language of \cite{AAEL}, this is the seed.

Next, we obtain density for the projective unitary group for the space spanned by $e_1, e_2, e_3$ using the ``bridge" $\sigma_3$.  
This relies on the ``Bridge Lemma" \cite[Lemma 8.2]{AAEL} and the ``Decoupling Lemma" \cite[Lemma 8.3]{AAEL}.
Continuing in this way, we obtain density for the whole space.
\end{proof}


\begin{remark}
Density results for braid group actions in the non-compact case have been considered by Kuperberg \cite{Kupdense} and Aharonov, Arad, Eban, and Landau \cite{AAEL}.  Both papers were partially motivated by quantum algorithms for the Potts model (see also \cite{Kupjones}), which is related to the Jones polynomial at various values of the quantum parameter.  When this parameter is not a root of unity, one is led to Hermitian vector spaces with mixed signatures and the corresponding Lie groups are then non-compact.

Kuperberg takes a more Lie theoretic approach, while Aharonov, Arad, Eban, and Landau develop some elementary linear algebra machinery.
%
%
Funar \cite{Funarmixeddense} obtains density results for actions of mapping class groups beyond braid groups. These actions often have image in indefinite unitary groups.
\end{remark}

 \bibliographystyle{plain}
\bibliography{bib_unitary}

\end{document}